\documentclass[10pt, amstex]{article}

\usepackage{amssymb}
\usepackage{fullpage}
\usepackage{amsthm}
\usepackage{amsmath}
\usepackage{mathrsfs}
\begin{document}
\newcommand{\beq}{\begin{eqnarray}}
\newcommand{\eeq}{\end{eqnarray}}
\newcommand{\beas}{\begin{eqnarray*}}
\newcommand{\enas}{\end{eqnarray*}}
\newcommand{\bea}{\begin{eqnarray}}
\newcommand{\ena}{\end{eqnarray}}
\newcommand{\D}{$\dagger$}
\newcommand{\dd}{$\ddagger$}
\newtheorem{theorem}{Theorem}[section]
\newtheorem{corollary}{Corollary}[section]
\newtheorem{conjecture}{Conjecture}[section]
\newtheorem{proposition}{Proposition}[section]
\newtheorem{lemma}{Lemma}[section]
\newtheorem{definition}{Definition}[section]
\newtheorem{condition}{Condition}[section]
\newtheorem{remark}{Remark}[section]
\newcommand{\pf}{\noindent {\bf Proof:} }
\title{$L^p$ bounds for a  central limit theorem with involutions}
%\runtitle{$L^p$ bounds for involutions}
\author{Subhankar Ghosh\footnote{Department of Mathematics, University of Southern California, Los Angeles, CA 90089, USA,\texttt{subhankg@usc.edu}}\\
\\
{\normalsize{\em University of Southern California } }}
\date{}
\maketitle

 \long\def\symbolfootnote[#1]#2{\begingroup%
\def\thefootnote{\fnsymbol{footnote}}\footnote[#1]{#2}\endgroup}
\symbolfootnote[0]{2000 {\em Mathematics Subject Classification}: Primary 60F25; Secondary 60F05,60C05.}
\symbolfootnote[0]{{\em Keywords}: CLT, Combinatorial Limit Theorems, Stein's method.}
\begin{abstract}
Let $E=((e_{ij}))_{n\times n}$ be a fixed array of real numbers such that $e_{ij}=e_{ji}, e_{ii}=0$ for $1\le i,j \le n$. Let the permutation group
be denoted by $S_n$ and the collection of involutions with no fixed points by $\Pi_n$, that is, $\Pi_n=\{\pi\in S_n: \pi^2=\mbox{id}, \pi(i)\neq i\,\forall i\}$ with $\mbox{id}$ denoting the identity permutation. For $\pi$ uniformly chosen from $\Pi_n$, let $Y_E=\sum_{i=1}^n e_{i\pi(i)}$ and $W=(Y_E-\mu_E)/\sigma_E$ where $\mu_E=E(Y_E)$ and $\sigma_E^2=\mbox{Var}(Y_E)$. Denoting by $F_W$ and $\Phi$ the distribution functions of $W$ and a $\mathcal{N}(0,1)$ variate respectively, we bound $||F_W-\Phi||_p$ for $ 1\le p\le \infty$ using Stein's method and the zero bias transformation. Optimal Berry-Esseen or $L^\infty$ bounds for the classical problem where $\pi$ is chosen uniformly from  $S_n$ were obtained by Bolthausen using Stein's method. Although in our case $\pi \in \Pi_n$ uniformly, the $L^p$ bounds we obtain are of similar form as Bolthausen's bound which holds for $p=\infty$. The difficulty in extending Bolthausen's method from $S_n$ to $\Pi_n$ arising due to the involution restriction is tackled by the use of zero bias transformations.
\end{abstract}
\section{Introduction}
Let $E=((e_{ij}))$ be an $n\times n$ array of real numbers. The study of combinatorial central limit theorems, that is, central limit theorems for random variables of the form
\beq \label{def-U}
Y_E=\sum_{i=1}^n e_{i\pi(i)} \eeq
focuses on the case where $\pi$ is a permutation chosen uniformly from a subset $A_n$ of $S_n$, the permutation group of order $n$. Some of the well studied choices for $A_n$ are $S_n$ itself \cite{hoeffding}\cite{ho}\cite{bolt}\cite{goldstein}, the collection $\Pi_n$ of fixed point free involutions \cite{goldstein2} \cite{schiff}, and the collection of permutations having one long cycle  \cite{kolchin}. The last two cases of $A_n$ are examples of distributions over $S_n$ that are constant on conjugacy classes, considered in \cite{goldstein-pattern}.

In this paper, we will be interested in the specific case of $A_n=\Pi_n$, where
\bea
\Pi_n=\{\pi\in S_n:\pi^2=\mbox{id}, \forall i:\pi(i)\not=i\}\label{def-pi-n}
\ena
with id denoting the identity permutation. However, before focusing on this choice, and the technicalities caused by restricting our permutations to be fixed point free involutions, we first briefly review the existing results pertaining to the much more studied case $A_n=S_n$, otherwise commonly known as the Hoeffding combinatorial CLT.

Approximating the distribution of $Y_E$ by the normal
when $\pi$ is chosen uniformly from $S_n$ began with the work of Wald and Wolfowitz \cite{ww}, who,
motivated by approximating null distributions for
permutation test statistics, proved
the central limit theorem as $n \rightarrow \infty$ for the case
where the factorization $e_{ij}=b_ic_j$ holds. In this special case, when $b_i$ are numerical characteristics
of some population, $k$ of the $c_j$'s are 1, and the remaining $n-k$ are zero, $Y_E$
has the distribution of a sum obtained by simple random sampling of size $k$ from the population. General arrays were handled in
the work of Hoeffding \cite{hoeffding}, and Motoo obtained Lindeberg-type conditions which are sufficient for the normal limit in \cite{motoo}.

A number of later authors refined these limiting results and obtained information on
the rate of convergence and bounds on the error in the normal approximation, typically
in the supremum or $L^\infty$ norm. Ho and Chen \cite{ho}
 and von Bahr \cite{von} derived $L^\infty$ bounds when the matrix $E$ is
random, the former using a concentration inequality approach and Stein's method, yielding the correct rate $O(n^{-1/2})$ under certain boundedness assumptions on $\sup_{i,j}|e_{i,j}|$. Goldstein \cite{goldstein-pattern}, employing the zero bias version of Stein's method obtained bounds of the correct order  with an explicit constant, but in terms of $\sup_{i,j}|e_{i,j}|$.
The work of Bolthausen \cite{bolt}, proceeding inductively, is the only one which yields an $L^\infty$ bound in terms of third moment type quantities on $E$ without the need for conditions on $\sup_{i,j}|e_{i,j}|$. More recently, Goldstein \cite{goldstein} obtained $L^1$ bounds for this case using zero biasing.
%Though an explicit constant was
%not obtained in \cite{bolt}, as we will see here for a related case, this inductive method can yield explicit constants when the %calculations are tracked carefully.

\par The case of $A_n=\Pi_n$ was considered much more recently. In \cite{schiff}, a permutation test is considered for a certain matched pair experiment designed to answer the question of whether there is an unusually high degree of similarity in a distinguished pairing, when there is some unknown background or baseline similarity between all pairs. In such a case, one considers $Y_E$ as in (\ref{def-U}) when $\pi$ is chosen uniformly from $\Pi_n$.

 Since the distribution of $Y_E$ is complicated, $L^\infty$ bounds for the error in normal approximation enables one to test the significance of the matched pair. The interested reader can look into \cite{goldstein2} for a discussion on similar applications.

 A bound to the normal for this case was provided in the $L^\infty$ norm by \cite{goldstein2}, with explicit constants along with the order, but under a boundedness assumption.  In this paper, we use techniques similar to \cite{bolt} and \cite{goldstein} to relax the conditions of \cite{goldstein2} so that $L^p$ bounds to the normal for the involution case can be obtained for possibly unbounded arrays $E$ also, in terms of third moment type quantities on the matrix $E$. In particular, in Theorem \ref{Lp-theorem} we show that
if $\pi\in\Pi_n$ uniformly, $W$ denotes the variable $Y_E$ appropriately standardized and $K_p$ is given explicitly by (\ref{def-Kp}), then the $L^p$ norm of the difference between $W$ and the normal satisfies
$$
||F_W-\Phi||_p \le K_p \frac{\beta_E}{n}
$$
for all $n \ge 9$, where $\beta_E$ is given in (\ref{def-beta}). This error bound yields a rate of $O(n^{-1/2})$ in the case of bounded arrays and as indicated in the appendix, for the bounded arrays this rate can not be improved uniformly over all arrays. Although the constant $K_p$ is quite large in magnitude to be applied for practical example, it is the first of its kind in the literature. Also we improve upon Goldstein and Rinott's \cite{goldstein2} result since we obtain bounds of order $O(n^{-1/2})$ under milder conditions like $\beta_E/\sqrt{n}$ being bounded instead of $\sup |e_{ij}|$ being bounded. It should be noted that the method applied here can be adopted to give $L^p$ estimates in Hoeffding's combinatorial CLT as well, and will yield a bound of the same from as the one obtained in \cite{bolt}.

While Bolthausen's method \cite{bolt} yields optimal results for the Hoeffding combinatorial CLT, that is the case of $A_n=S_n$, it is not immediately clear to the author if it can be extended to other classes of permutations including the case of involutions that is $A_n=\Pi_n$. The main problem is that the auxiliary permutations $\pi_1,\pi_2,\pi_3$ considered in page 383 of \cite{bolt} are unrestricted, but similar permutations for the involutions CLT have to be involutions without a fixed point which makes the construction harder. As we shall see, this difficulty can be overcome in a natural way by using the zero biasing. The auxiliary variables produced by Proposition \ref{prop:uniform-interpolation} are absolutely continuous with respect to $Y_E$ as in (\ref{def-U}) with $\pi\in\Pi_n$. This, in particular, ensures that the corresponding auxiliary permutations we obtain are involutions. A second novelty of using the zero bias transformation is that it yields not only the optimal $L^\infty$ or Berry-Esseen bounds, but general $L^p$ bounds holding for all $1\le p\le\infty$.

Another work where $\pi$ is not uniform is that of \cite{kolchin}, where the permutations are uniform over those having one long cycle, and \cite{goldstein-pattern}, where $L^\infty$ bounds are derived under a boundedness condition for the case of a permutation distribution constant on conjugacy classes having no fixed points, which generalizes both the involution and long cycle cases.

The paper is organized as follows. In Section 2, we introduce some notation and state our main result. In Section 3 the basic idea of the zero bias transformation is reviewed, and an outline is provided that illustrates how to obtain zero bias couplings in some cases of interest. In Section 4, $L^1$ bounds are obtained. Lastly, in Section 5 we use the calculations in Section 4 along with the recursive argument of \cite{bolt} to obtain $L^\infty$ bounds. From the $L^1$ and $L^\infty$ bounds, the following simple inequality allows for the computation of $L^p$ bounds for all intermediate $p \in (1,\infty)$,
\beq
||f||^p_p\le ||f||_\infty^{p-1}||f||_1.\label{ineq-p}
\eeq

\section{Notation and statement of main result}
For $n$ an even positive integer, let $\pi$ be a permutation
chosen uniformly from $\Pi_n$ in (\ref{def-pi-n}), the set of involutions with no fixed points.
Since for $\pi \in \Pi_n$ the terms $e_{i\pi(i)}$ and $e_{\pi(i)i}$
always appear together in (\ref{def-U}), and $e_{ii}$ never appears, we may assume
without loss of generality that
\bea \label{e-symm-ii}
e_{ij}=e_{ji} \quad \mbox{and that} \quad e_{ii}=0 \quad \mbox{for
all $i,j=1,2,\ldots,n$.} \ena

For an array $E=((e_{ij}))_{1\le i,j \le n}$ satisfying the conditions in (\ref{e-symm-ii}), define
$$
e_{i + }=\sum_{j=1}^n e_{ij}, \quad e_{+ j}=\sum_{i=1}^n e_{ij}\quad \mbox{and }
e_{++}=\sum_{i,j=1}^n e_{ij}.
$$
Then, as shown in \cite{goldstein2},
\beq
\label{mean-var-Ye}
\begin{array}{ccl}
\mu_E&=&EY_E=\frac{e_{++}}{n-1}\\
\sigma_E^2&=&\frac{2}{(n-1)(n-3)}\left((n-2)\sum_{1\le i,j \le
n}e_{ij}^2 + \frac{1}{n-1}e_{++}^2-2\sum_{i=1}^n e_{i+}^2 \right).
\end{array}
\eeq
Again following \cite{goldstein2}, letting
\bea
\label{def-ghat}
{\widehat e}_{ij}=\left\{
\begin{array}{cc}
e_{ij} -\frac{e_{i +}}{n-2} - \frac{e_{+j}}{n-2} + \frac{e_{+
+}}{(n-1)(n-2)}&  i \not = j\\
0 & i=j,
\end{array}
\right.
\ena
 we have
 \bea \label{hat-e-sums}
\widehat{e}_{+ i}= \widehat{e}_{j+}=\widehat{e}_{++}=0 \quad \mbox{for all
$i,j=1,\ldots,n$.}
\ena
and that
\beas  Y_{\widehat{E}}=Y_E-\mu_E \enas
where $\widehat{E}$ is obtained from $E$ by (\ref{def-ghat}).

We consider bounds to the normal ${\cal N}(0,1)$ for the
standardized variable \bea \label{def-Y} W=\frac{Y_E -
\mu_E}{\sigma_E}. \ena
Since $Y_E$ and $Y_{\widehat{E}}$ differ by a constant,
(\ref{mean-var-Ye}) and (\ref{hat-e-sums}) yield
\bea \label{sigmaD}
\sigma_E^2=\sigma_{\widehat{E}}^2 = \frac{2(n-2) }{(n-1)(n-3)}
\sum_{1\le i,j \le n}\widehat{e}_{ij}^2.
\ena
In particular, the mean zero, variance 1 random variable $W$
in (\ref{def-Y}) can be written as \beq \label{def-Yd} W=\sum_{i=1}^n
d_{i\pi(i)} \quad \mbox{where $d_{ij}=\widehat{e}_{ij}/\sigma_E$}, \eeq and
moreover, the array $\widehat{E}$ inherits properties (\ref{e-symm-ii}) and (\ref{hat-e-sums}) from $E$, as then does $D$ from $\widehat{E}$.
For any array $E=((e_{ij}))_{n\times n}$, let
\beq\label{def-beta}
\beta_E=\sum_{i\neq j} \frac{|\widehat{e}_{ij}|^3}{\sigma_E^3}\quad\mbox{with $\widehat{e}_{ij}$ as in (\ref{def-ghat}) and $\sigma_E^2$ as in (\ref{mean-var-Ye}).}\eeq  In order to guarantee
that $W$ is a well defined random variable, that is, that $\sigma_E^2>0$, we henceforth impose the following
condition without further mention.

\begin{condition}\label{rem-cond-e}
For some $i \not =j$, the value $\widehat{e}_{ij}\not=0$, or, equivalently,
$$
e_{ij} -\frac{e_{i +}}{n-2} - \frac{e_{+j}}{n-2} + \frac{e_{+
+}}{(n-1)(n-2)}\not=0\quad\mbox{for some $i\not=j$.}
$$
\end{condition}
\noindent Since $d_{ij}$ and $\widehat{e}_{ij}$ are linearly related, Condition \ref{rem-cond-e} is equivalent
to the condition that $d_{ij} \not =0$ for some $i \not =j$.

We provide bounds on the accuracy of the normal approximation of $W$ using the zero bias transformation,
introduced in \cite{goldstein1}. For any random variable $W$ with mean zero and variance $\sigma^2$, there exists a unique distribution for a random variable $W^*$ with the property that for any
differentiable function $f$
\beq
E[Wf(W)]=\sigma^2E f'(W^*).\label{def-zbias}
\eeq
The distribution of the random variate $W^*$ is called the zero bias transform of the distribution of $W$.
 From Stein's original lemma \cite{stein}, $W$ is ${\cal N}(0,\sigma^2)$ if
and only if $W^*$ is ${\cal N}(0,\sigma^2)$, that is, the normal distribution is the unique fixed point
  of the zero bias transformation. This gives rise to the intuition that if $W$ and $W^*$ are close, in some appropriate sense, then the distribution of $W$ should be close to the normal. That this intuition is indeed true can be seen, for instance, in the following result from \cite{goldstein}: if $W$ and $W^*$ are on a joint space, with $W^*$ having the $W$ zero bias distribution, then
\beq
||F_W-\Phi||_1 \le 2E|W^*-W|.\label{l1result}
\eeq
In (\ref{l1result}), $F_W$ and $\Phi$ denote the distribution functions of $W$ and that of a standard normal variate respectively, and $||\cdot||_p$ denotes the $L^p$ norm. We call a construction of $W$ and $W^*$ on a joint space a zero bias coupling of $W$ to $W^*$.

The following is our main result which we prove using zero bias coupling.
\begin{theorem}
\label{Lp-theorem}  Let $E=((e_{ij}))_{n\times n}$ be an array satisfying $e_{ij}=e_{ji},e_{ii}=0\,\forall i, j$, and let $\pi$ be
an involution chosen uniformly from $\Pi_n$. If
$$
Y_E=\sum_{i=1}^n e_{i \pi(i)},
$$
and $W=(Y_E-\mu_E)/\sigma_E$, then for $n \ge 9$ and $p \in [1,\infty]$, with $\beta_E$ as in (\ref{def-beta}), we have
$$
||F_W-\Phi||_p \le K_p \frac{\beta_E}{n}.
$$
Here $F_W$ denotes the distribution function of $W$, $\Phi$ is the distribution function of a standard normal variate and
\bea
\label{def-Kp}
K_p=(379)^{1/p}(61,702,446)^{1-1/p}.
\ena
\end{theorem}
As $W$ in the theorem is given by (\ref{def-Yd}) with
\bea \label{D-conditions}
d_{ij}=d_{ji}, d_{ii}=0, d_{i+}=0 \quad \mbox{and} \quad \sigma_D^2=1;\,\beta_E=\beta_D
\ena
we assume in what follows that all subsequent occurrences of $((d_{ij}))$
satisfy these conditions and instead of working with $E$ work with the centered and scaled array $D$ only.

In the next section, we review construction of zero bias couplings in certain cases of interest including the present problem.

\section{Zero bias transformation}
\par We prove Theorem \ref{Lp-theorem} by constructing a zero bias
coupling using a Stein pair, that is, a pair of random variables
$(W,W')$ which are exchangeable and satisfy \beq
\label{linearity-condition-lambda} E(W-W^\prime|W)=\lambda W \quad \mbox{for some $\lambda \in (0,1)$}; \eeq
 see \cite{chen} for more on Stein pairs.
\par As shown in \cite{goldstein1}, for any mean zero, variance $\sigma^2$ random variable $W$,
there exists a distribution for a random variable $W^*$ satisfying (\ref{def-zbias}). Nevertheless, constructing useful couplings of $W$ and $W^*$ for particular examples may be difficult. In some cases, however, as in ours,
the following proposition from \cite{goldstein} may be applied.

\begin{proposition}
\label{prop:uniform-interpolation}
Let $W,W^\prime$ be an exchangeable pair with
$Var(W)=\sigma^2\in(0,\infty)$ and distribution $F(w,w^\prime)$
satisfying the linearity condition
(\ref{linearity-condition-lambda}). Then
\bea \label{2lambdasigma2}
E(W-W^\prime)^2=2\lambda\sigma^2,
\ena
and when $(W^\dag,W^\ddag)$ has the joint distribution \beq
dF^\dag(w,w^\prime)=\frac{(w-w^\prime)^2}{E(W-W^\prime)^2}dF(w,w^\prime)\label{def-sqbias}
\eeq and $U\sim\mathcal{U}$[0,1] is independent of
$W^\dag,W^\ddag$, the variable
\beas
W^*=UW^\dag+(1-U)W^\ddag
\enas
has the $W$-zero biased distribution.
\end{proposition}
\par

In particular, construction of zero bias couplings are typically possible
when Stein pairs exist. We review the construction of
$(W^\dag,W^\ddag)$ from $(W,W^\prime)$ as outlined in
\cite{goldstein}.
Suppose we have a Stein pair $(W,W^\prime)$
which is a function of some collection of random variables
$\{\Xi_\alpha,\alpha\in\mathbf{\chi}\}$, and that for a possibly
random index set ${\bf I} \subset \mathbf{\chi}$, independent of
$\{\Xi_\alpha,\alpha\in \mathbf{\chi}\}$, the difference
$W-W^\prime$ depends only on ${\bf I}$ and on $\{\Xi_\alpha,\alpha\in \chi_{\bf {\footnotesize I}}\}$, where $\chi_\mathbf{I}\subset\chi$ depends on $\mathbf{I}$. That
is, for some function $b(\mathbf{i}, \Xi_\alpha,\alpha\in\chi_{\bf i})$ defined on ${\bf i} \subset \mathbf{\chi}$, and ${\bf I}$
a random index set,
\beq W-W^\prime=b({\bf I},\Xi_\alpha,\alpha\in\chi_{\bf
I}).\label{bdef} \eeq
\par We now show how, under this framework, the pair $(W,W^\prime)$ can
be constructed; the pair $(W^\dagger,W^\ddagger)$ will then be constructed in
a similar fashion. First generate $\mathbf{I}$,
then independently generate
$\{\Xi_\alpha,\alpha\in\chi_{\mathbf I}\}$ and finally
$\{\Xi_\alpha,\alpha\in \chi_\mathbf{I}^c\}$ conditioned on
$\{\Xi_\alpha,\alpha\in\chi_\mathbf{I}\}$. That is, first generate the indices ${\bf I}$ on which the difference $W-W'$ depends, then the underlying variables $\Xi_\alpha, \alpha \in \chi_{\bf I}$
which make up that difference, and lastly the remaining variables.
This construction corresponds to the following factorization of the joint
distribution of ${\bf I}$ and $\{\Xi_\alpha, \alpha \in \chi\}$ as the
product \beq
\,dF(\mathbf{i},\xi_\alpha,\alpha\in\chi)=P(\mathbf{I=i})\,
dF_{\mathbf{i}}(\xi_\alpha,\alpha\in\chi_\mathbf{i})
\,dF_{\mathbf{i}^c|\mathbf{i}}(\xi_\alpha,\alpha\notin\chi_\mathbf{i}|\xi_\alpha,\alpha\in\chi_\mathbf{i}).
\label{factor1} \eeq

For $dF^\dagger$ we consider the joint distribution of ${\bf
I}$ and $\{\Xi_\alpha, \alpha \in \chi\}$, biased by the squared difference
$(w-w^\prime)^2$, is \beq
dF^\dagger(\mathbf{i},\xi_\alpha,\alpha\in\chi)=
\frac{(w-w^\prime)^2}{E(W-W^\prime)^2}dF(\mathbf{i},\xi_\alpha,\alpha\in\chi).
\label{sqbias1} \eeq
From  (\ref{2lambdasigma2}), (\ref{bdef}) and the independence of ${\bf I}$
and $\{\Xi_\alpha, \alpha \in \chi\}$ we obtain \beq
\sum_{\mathbf{i}\subset\chi}
P(\mathbf{I=i})Eb^2(\mathbf{i},\Xi_\alpha,\alpha\in\chi_\mathbf{i})=2\lambda\sigma^2.
\eeq Hence we can define a probability distribution for an index set
$\mathbf{I^\dag}$ by
\bea \label{def-ri}
P(\mathbf{I^\dag=i})= \frac{r_\mathbf{i}}{2\lambda\sigma^2} \quad
\mbox{where} \quad
r_{\mathbf{i}}=P(\mathbf{I=i})Eb^2(\mathbf{i},\Xi_\alpha,\alpha\in\chi_\mathbf{i}).
\ena
From (\ref{bdef}), (\ref{sqbias1}) and (\ref{def-ri}), we obtain \beq
dF^\dagger(\mathbf{i},\xi_\alpha,\alpha\in\chi)=\frac{b^2(\mathbf{i},\xi_\alpha,\alpha\in\chi_\mathbf{i})}{2\lambda\sigma^2}P(\mathbf{I=i})dF_{\mathbf{i}}(\xi_\alpha,\alpha\in\chi_\mathbf{i})dF_{\mathbf{i}^c|\mathbf{i}}(\xi_\alpha,\alpha\notin\chi_\mathbf{i}|\xi_\alpha,\alpha\in\chi_\mathbf{i})\nonumber\\
=\frac{r_\mathbf{i}}{2\lambda\sigma^2}\frac{b^2(\mathbf{i},\xi_\alpha,\alpha\in\chi_\mathbf{i})}{Eb^2(\mathbf{i},\Xi_\alpha,\alpha\in\chi_\mathbf{i})} dF_\mathbf{i}(\xi_\alpha,\alpha\in\chi_\mathbf{i})dF_{\mathbf{i}^c|\mathbf{i}}(\xi_\alpha,\alpha\notin\chi_\mathbf{i}|\xi_\alpha,\alpha\in\chi_\mathbf{i})\nonumber\\
=P(\mathbf{I^\dag=i})dF^\dag_\mathbf{i}(\xi_\alpha,\alpha\in\chi_\mathbf{i})dF_{\mathbf{i}^c|\mathbf{i}}(\xi_\alpha,\alpha\notin\chi_\mathbf{i}|\xi_\alpha,\alpha\in\chi_\mathbf{i})\label{factor2}
\eeq where
\bea \label{def-dFi}
dF_\mathbf{i}^\dag(\xi_\alpha,\alpha\in\chi_\mathbf{i})=\frac{b^2(\mathbf{i},\xi_\alpha,\alpha\in\chi_\mathbf{i})}{Eb^2(\mathbf{i},\Xi_\alpha,\alpha\in\chi_\mathbf{i})}dF_\mathbf{i}(\xi_\alpha,\alpha\in\chi_\mathbf{i}).
\ena
\par Note that (\ref{factor2}) gives a representation of the distribution $F^\dag$
which is of the same form as (\ref{factor1}). The parallel forms
of $F$ and $F^\dagger$ allow us
to generate variables having distribution $F^\dag$ parallel to that for distribution $F$;
first generate the random index $\mathbf{I}^\dag$, then
$\{\Xi_\alpha^\dag,\alpha\in\chi_\mathbf{I}^\dag\}$ according to
$dF_{\mathbf{I}}^\dag$, and lastly the remaining variables
according to $dF_{\mathbf{I}^c|\mathbf{I}}(\xi_\alpha,\alpha\notin\chi_\mathbf{i}|\xi_\alpha,\alpha\in\chi_\mathbf{i})$.

That the last step is the same in the construction of pairs of variables having the $F$ and $F^\dagger$ distribution
allows an opportunity for a coupling between $(W,W')$ and $(W^\dagger, W^\ddagger)$ to be achieved; from Proposition \ref{prop:uniform-interpolation}, this suffices to couple $W$ and $W^*$. One coupling may be accomplished using the following outline, as was done in \cite{goldstein}.
First generate ${\bf I}$ and $\{\Xi_\alpha, \alpha \in \chi\}$, yielding the pair $W,W'$.
Next generate ${\bf I}^\dagger$ and then the variables $\{\Xi_\alpha^\dagger:\alpha\in\chi_{\mathbf{I}^\dagger}\}$ following $dF^\dag_\mathbf{i}$ if the realization of $\mathbf{I}^\dag$ is $\mathbf{i}$.
Lastly, when constructing the remaining variables that is $\{\Xi_\alpha^\dagger: \alpha \notin \chi_{\mathbf{i}}\}$ which
make up $W^\dagger, W^\ddagger$, use
as much of the previously generated variables $\{\Xi_\alpha, \alpha \notin \chi_\mathbf{i}\}$ as possible so that the pairs
$(W,W')$ and $(W^\dagger,W^\ddagger)$ will be close.

\par We review the construction of $(W,W')$ in \cite{goldstein2}
for the case at hand, and then show
how it agrees with the outline above. For distinct $i,j  \in
\{1,\ldots,n\}$ let $\tau_{i,j}$ be the permutation which
transposes the elements $i$ and $j$, that is, $\tau_{i,j}(i)=j,
\tau_{i,j}(j)=i$ and $\tau_{i,j}(k)=k$ for all $k \not \in
\{i,j,\}$. Furthermore, given $\pi \in \Pi_n$, let \beas
 \alpha^\pi_{i,j}=\tau_{i,\pi(j)} \tau_{j,\pi(i)}.
\enas Note that for any given $\pi \in \Pi_n$, the permutation $\pi
\alpha^\pi_{i,j}$ will again belong to $\Pi_n$. In particular, whereas
$\pi$ has the cycle(s) $(i,\pi(i))$ and $(j,\pi(j))$, the
permutation $\pi \alpha^\pi_{i,j}$ has cycle(s) $(i,j)$ and
$(\pi(i),\pi(j))$, and all other cycles in common with $\pi$. Now,
with $\pi$ chosen uniformly from $\Pi_n$ and $W$ given by
(\ref{def-Yd}), we construct an exchangeable pair $(W,W')$, as
in \cite{goldstein2}, as follows. Choose two distinct indices
$\mathbf{I}=(I,J)$ uniformly from $\{1,2,\ldots,n\}$, that is, having distribution
\bea \label{dist-I,J}
P({\bf I}={\bf i})=P(I=i,J=j)=\frac{1}{n(n-1)}{\bf 1}(i \not = j)
\ena
and let \beq
\label{pi-prime-from-pi} \pi^\prime =  \pi \alpha^\pi_{I,J}. \eeq
Since $(I,J)$ is chosen uniformly over all distinct pairs, $\pi$ and
$\pi^\prime$ are exchangeable, and hence, letting $W^\prime=\sum
d_{i\pi^\prime(i)}$, so are $(W,W^\prime)$. Moreover, as $\pi^\prime$
has the cycle(s) $(I,J)$ and $(\pi(I),\pi(J))$, and shares all other
cycles with $\pi$, we have \beq
W-W^\prime=2(d_{I\pi(I)}+d_{J\pi(J)}-(d_{IJ}+d_{\pi(I)\pi(J)})).\label{bref}
\eeq
Using (\ref{bref}) it is shown in \cite{goldstein2} that
\bea\label{lambda-inv}
E(W-W^\prime|W)=\frac{4}{n}W,
\ena
that is, (\ref{linearity-condition-lambda}) is satisfied with
$\lambda=4/n$.

To put this construction in the framework above, so to be able to apply the decomposition (\ref{factor2}) for the construction of a zero bias coupling,
let $\chi=\{1,2,\ldots,n\}$ and
$\Xi_\alpha=\pi(\alpha)$. From (\ref{bref}), we see that the difference $W-W'$ depends on
a pair of randomly chosen indices, and their images. Hence, regarding these indices,
let ${\bf i} =(i,j) \in \chi^2$ and let ${\bf I}=(I,J)$ where $I$ and
$J$ have joint distribution given in (\ref{dist-I,J}), specifying $P({\bf I}={\bf i})$, the
first term in (\ref{factor1}). Also, in this case $\chi_\mathbf{i}=\mathbf{i}$. Next, for given
distinct $i,j$ and $\pi \in \Pi_n$, we have $\pi(i) \not = \pi(j), \pi(i) \not =i$ and $\pi(j) \not = j$. As $\pi$ is chosen uniformly from $\Pi_n$, the distribution of the images $k=\xi_i$ and
$l=\xi_j$ of distinct $i$ and $j$ under $\pi$ is given by
\beq
dF_{{\bf i}}(\xi_\alpha, \alpha \in {\bf i})=dF_{i,j}(k,l) \propto {\bf 1}(k \not = l, k \not = i, l \not = j)  \quad
\mbox{for $i \not = j$},\label{def-df}
\eeq
specifying the second term of (\ref{factor1}). The last term in (\ref{factor1}) is given by
\bea
dF_{{\bf i}^c|{\bf i}}(\xi_\alpha, \alpha \not \in {\bf i}|\xi_\alpha, \alpha \in {\bf i})&=&\frac{P(\pi(i)=k, \pi(j)=l,\pi(\alpha)=\xi_\alpha, \alpha \not \in \{i,j,k,l\})}{P(\pi(i)=k, \pi(j)=l)},\label{cond-general}
\ena
when $i\neq j$ and $k\neq l$.
 The equality in (\ref{cond-general}) follows from the fact that $\pi$ is an involution, implying $\xi_k=\pi(k)=\pi(\pi(i))=i$ and similarly for $\xi_l=j$, and thus
\beas
&&\{\pi(i)=k, \pi(j)=l,\pi(\alpha)=\xi_\alpha, \alpha \not \in \{i,j\}\}\\
&=&\{\pi(i)=k, \pi(j)=l,\pi(k)=i,\pi(l)=j,\pi(\alpha)=\xi_\alpha,\alpha \not \in \{i,j,k,l\}\}\\
&=&\{\pi(i)=k, \pi(j)=l,\pi(\alpha)=\xi_\alpha, \alpha \not \in \{i,j,k,l\}\}.
\enas
The conditional distribution in (\ref{cond-general}) can be simplified further by considering the two cases $k=j$ or equivalently $l=i$ and $k\neq j$ or equivalently $|\{i,j,k,l\}|=4$ separately. If $k=j$ then we have $l=i$ and thus we obtain
\bea
dF_{{\bf i}^c|{\bf i}}(\xi_\alpha, \alpha \not \in {\bf i}|\xi_i=j,\xi_j=i)&=&\frac{P(\pi(i)=j, \pi(j)=i,\pi(\alpha)=\xi_\alpha, \alpha \not \in \{i,j\})}{P(\pi(i)=j, \pi(j)=i)}\nonumber\\
&=&\frac{P(\pi(i)=j,\pi(\alpha)=\xi_\alpha, \alpha \not \in \{i,j\})}{P(\pi(i)=j)}\nonumber\\
&=&\frac{|\Pi_n|^{-1}}{(n-1)^{-1}}=|\Pi_{n-2}|^{-1}.\label{cond-uniform-2}
\ena
For (\ref{cond-uniform-2}), we note $P(\pi(i)=j)=1/(n-1)$ since $\pi$ is chosen uniformly from $\Pi_n$ .
 The last equality in (\ref{cond-uniform-2}) simply indicates that once we fix $\pi(i)=j$, that is the cycle $(i,j)$ in $\pi$ we have to choose an involution uniformly at random from the rest of the indices that is $\Pi_{n-2}$ to obtain $\pi$ in its entirety. This argument yields the recursion which we use in (\ref{cond-uniform}) later
 \beas
 |\Pi_n|=(n-1)|\Pi_{n-2}|.
 \enas
 When we have $k\neq j$ and hence $l\neq  i$ or equivalently $|\{i,j,k,l\}|=4$ in (\ref{cond-general}),  we obtain
\bea
dF_{{\bf i}^c|{\bf i}}(\xi_\alpha, \alpha \not \in {\bf i}|\xi_\alpha, \alpha \in {\bf i})&=&\frac{P(\pi(i)=k, \pi(j)=l,\pi(\alpha)=\xi_\alpha, \alpha \not \in \{i,j,k,l\})}{P(\pi(i)=k, \pi(j)=l)}\nonumber\\
&=& \frac{|\Pi_n|^{-1}}{((n-1)(n-3))^{-1}}=|\Pi_{n-4}|^{-1}.\label{cond-uniform}
\ena

In (\ref{cond-uniform}) we have used the following equality which follows from the fact that $\pi$ is an involution chosen uniformly at random.
 $$
 P(\pi(i)=k,\pi(j)=l)=P(\pi(i)=k)P(\pi(j)=l|\pi(i)=k)=\frac{1}{(n-3)(n-1)}.
 $$
 From (\ref{cond-uniform}) we see that the conditional distribution $dF_{{\bf i}^c|{\bf i}}$ is uniform over all values of $\xi_\alpha$ for $\alpha \in \chi$ for which $\xi_i=k, \xi_j=l$ and $P(\pi(\alpha)=\xi_\alpha, \alpha \in \chi)>0$ when $|\{i,j,k,l\}|=4$.  Hence we may construct $(W,W')$ following (\ref{factor1}), that is, first choosing ${\bf I}=(I,J)$,
then the images ${\bf K}=(K,L)$ of $I$ and $J$ under $\pi$ and $\pi'$, then the remaining images uniformly over all possible values for which the resulting permutations lie in $\Pi_n$.

We may construct $(W^\dagger, W^\ddagger)$ from (\ref{factor2}) quite easily now. In view of (\ref{bdef}) and (\ref{bref}),
for pairs of distinct indices $i,j$, let
\beq
b(i,j,\pi_i, \pi_j)=2(d_{i\pi_i}+
d_{j\pi_j}-(d_{ij}+d_{\pi_i\pi_j})),\label{bdef2} \eeq
where again we may also let $k=\xi_i$ and $l=\xi_j$. Now, considering the first two factors in (\ref{factor2}), using (\ref{def-ri}), (\ref{def-dFi}), (\ref{dist-I,J}) and (\ref{def-df}) we obtain,
\bea
\nonumber P(\mathbf{I^\dag=i})dF^\dag_\mathbf{i}(\xi_\alpha,\alpha\in\mathbf{i}) &=& P(I^\dagger=i,J^\dagger=j)
dF^\dag_{i,j}(K^\dagger=k,L^\dagger=l)\\
\nonumber &\propto& P(I=i,J=j)[d_{ik}+d_{jl}-(d_{ij}+d_{kl})]^2  {\bf 1}(k \not = l, k \not = i, l \not = j)\\
&\propto& [d_{ik}+d_{jl}-(d_{ij}+d_{kl})]^2  {\bf 1}(k \not = l, k \not = i, l \not = j, i \not = j).
\label{factor2-firsttwo}\ena
By Lemma \ref{lemma-cn}, below, relation (\ref{factor2-firsttwo}) specifies a joint distribution, say $p(i,j,k,l)$, on the pairs $\mathbf{I^\dag}=(I^\dag, J^\dag)$ and
their images $(\xi_{I^\dag},\xi_{J^\dag})=(K^\dag,L^\dag):=\mathbf{K}^\dag$, say. Since when $j =k$, or equivalently $i=l$ the square term vanishes, we may write
\beq
p(i,j,k,l)=c_n[d_{ik}+d_{jl}-(d_{ij}+d_{kl})]^2{\bf 1}(|\{i,j,k,l\}|=4),\label{def-p}
\eeq
where the constant of proportionality $c_n$ is provided in Lemma \ref{lemma-cn}. Next, note that since $i,j,k,l$ have to be distinct in the definition (\ref{def-p}), the third
term in (\ref{factor2}), $dF_{\mathbf{i}^c|\mathbf{i}}(\xi_\alpha,\alpha\notin\mathbf{i}|\xi_\alpha,\alpha\in\mathbf{i})$, reduces to uniform distribution over all values of $\xi_\alpha$ for $\alpha \in \chi$ for which $\xi_i=k, \xi_j=l$ and $P(\pi(\alpha)=\xi_\alpha, \alpha \in \chi)>0$ owing to (\ref{cond-uniform}). Once we form $(W^\dag,W^\ddag)$ following the square bias distribution (\ref{def-sqbias}), its easy to produce $W^*$ which has the $W$ zero bias distribution using Proposition \ref{prop:uniform-interpolation}. We summarize the conclusions above in the following lemma.

\begin{lemma}
\label{daggers-construction-lemma}
Let $dF(w,w')$ be the joint distribution of a Stein pair $(W,W')$ where
$$
W=\sum_{i=1}^n d_{i \pi(i)} \quad \mbox{and} \quad W'=\sum_{i=1}^n d_{i \pi'(i)}
$$
with $\pi$ chosen uniformly from $\Pi_n$ and $\pi'$ as in (\ref{pi-prime-from-pi})
with $I,J$ having distribution (\ref{dist-I,J}). Then
a pair $(W^\dagger,W^\ddagger)$ with the square bias distribution (\ref{def-sqbias}) can be constructed by setting
$$
W^\dagger=\sum_{i=1}^n d_{i \pi^\dagger(i)} \quad \mbox{and} \quad W^\ddagger=\sum_{i=1}^n d_{i \pi^\ddagger(i)}
$$
where $\pi^\dagger$ and $\pi^\ddagger=\pi^\dagger\alpha^{\pi^\dag}_{I^\dagger,J^\dagger}$ are constructed by first sampling ${\bf I}^\dagger = (I^\dagger,J^\dagger)$ and the respective images ${\bf K}^\dagger=(K^\dagger,L^\dagger)$ under $\pi^\dagger$ according to (\ref{def-p}) and then selecting the remaining images of $\pi^\dagger$ uniformly from among the choices for which it lies in $\Pi_n$. Furthermore if $U\sim\mathcal{U}[0,1]$ is independent of $(W^\dag,W^\ddag)$, then $W^*=UW^\dag+(1-U)W^\ddag$ has the $W$ zero bias distribution. Hence, if $\pi$ and
$\pi^\dag$ are constructed on a common space, then so are $W$ and $W^*$.
\end{lemma}

Given the permutations $\pi$ and $\pi'$ from which the pair $(W,W')$ is constructed, we would like to form
 $(W^\dag,W^\ddag)$ as close to $(W,W^\prime)$ as possible, thus making $W$ close to $W^*$. Towards this end, we follow
 the construction noted after (\ref{factor2}), using many of the already chosen variables which form $\pi$ and $\pi'$
 to make the two pairs close. In particular, begin the construction of $(W^\dag,W^\ddag)$ by choosing ${\bf I}^\dagger=(I^\dagger,J^\dagger)$ and ${\bf K}^\dagger=(K^\dagger, L^\dagger)$ with joint distribution (\ref{def-p}), independent of $\pi$ and $\pi'$.
 Let $R_1=|\{\pi(I^\dag),\pi(J^\dag)\}\cap\{K^\dag,L^\dag\}|$ and $R_2=|\{\pi(I^\dag),\pi(K^\dag)\}\cap\{J^\dag,L^\dag\}|$; clearly $R_1,R_2\in\{0,1,2\}$. Define $\pi^\dag$ by
 \beq\label{def-pi-dag}
\pi^\dagger=\left\{\begin{array}{clcr}
\pi\alpha^\pi_{J^\dag,L^\dag}&\mbox{if $\pi(I^\dag)=K^\dag$ and $\pi(J^\dag)\not=L^\dag$ }&\mbox{hence $(R_1,R_2)=(1,0)$}\\
\pi\alpha^\pi_{I^\dag,K^\dag}&\mbox{if $\pi(I^\dag)\not= K^\dag$ and $\pi(J^\dag)=L^\dag$}&\mbox{hence $(R_1,R_2)=(1,0)$}\\
\pi\alpha^\pi_{J^\dag,K^\dag}\tau_{I^\dag,J^\dag}\tau_{K^\dag,L^\dag}& \mbox{if $\pi(I^\dag)=L^\dag$ and $\pi(J^\dag)\not=K^\dag$}&\mbox{hence $(R_1,R_2)=(1,1)$}\\
\pi\alpha^\pi_{I^\dag,L^\dag}\tau_{I^\dag,J^\dag}\tau_{K^\dag,L^\dag}& \mbox{if $\pi(I^\dag)\not =L^\dag$ and $\pi(J^\dag)=K^\dag$}&\mbox{hence $(R_1,R_2)=(1,1)$}\\
\pi\alpha^\pi_{K^\dag,L^\dag}\tau_{I^\dag,L^\dag}\tau_{J^\dag,K^\dag}&\mbox{if $\pi(I^\dag)=J^\dag$ and $\pi(K^\dag)\not=L^\dag$}&\mbox{hence $(R_1,R_2)=(0,1)$}\\
\pi\alpha^\pi_{I^\dag,J^\dag}\tau_{I^\dag,L^\dag}\tau_{J^\dag,K^\dag}&\mbox{if $\pi(I^\dag)\not=J^\dag$ and $\pi(K^\dag)=L^\dag$}&\mbox{hence $(R_1,R_2)=(0,1)$}\\
\pi&\mbox{if $\pi(I^\dag)=K^\dag$ and $\pi(J^\dag)=L^\dag$}&\mbox{hence $(R_1,R_2)=(2,0)$}\\
\pi\tau_{I^\dag,L^\dag}\tau_{J^\dag,K^\dag}&\mbox{if $\pi(I^\dag)=J^\dag$ and $\pi(K^\dag)=L^\dag$}&\mbox{hence $(R_1,R_2)=(0,2)$}\\
\pi\tau_{I^\dag,J^\dag}\tau_{K^\dag,L^\dag}&\mbox{if $\pi(I^\dag)=L^\dag$ and $\pi(J^\dag)=K^\dag$}&\mbox{hence $(R_1,R_2)=(2,2)$}\\
\pi\alpha^\pi_{I^\dag,K^\dag}\alpha^\pi_{J^\dag,L^\dag}&\mbox{ when $R_1=R_2=0$}.
\end{array}\right.\eeq
 The partition in display (\ref{def-pi-dag}) is based on the possible values of $(R_1,R_2)$; it does not include the cases where $(R_1,R_2)=(2,1)$ or $(R_1,R_2)=(1,2)$ because these two events are impossible. If $R_1=2$ and $\pi(I^\dag)=K^\dag,\pi(J^\dag)=L^\dag$ then $R_2=0$ while if $\pi(I^\dag)=L^\dag,\pi(J^\dag)=K^\dag$ then $R_2=2$. Similarly one can rule out $(R_1,R_2)=(1,2)$. Similar arguments show us that the cases described above are indeed exhaustive.
\par Clearly any two cases in (\ref{def-pi-dag}) with differing values of $(R_1,R_2)$ tuple are exclusive. Also, any two cases with the same tuple value, e.g. cases one and two, are also exclusive. For example, in case one we have $\pi(I^\dag)=K^\dag$ whereas in case two we have $\pi(I^\dag)\not=K^\dag$ making these two cases disjoint. In summary $\pi^\dag$ is well defined, and
this construction specifies the pairs $(\pi,\pi')$ and $(\pi^\dagger,\pi^\ddagger)$ on the same space. The following lemma shows that the $\pi^\dag$ so obtained is an involution.
\begin{lemma}\label{lemma-pi-dag-involution}
For $\pi\in\Pi_n$, the permutation $\pi^\dag$ defined in (\ref{def-pi-dag}) belongs to $\Pi_n$ and has the cycles $(I^\dag,K^\dag)$ and $(J^\dag,L^\dag)$. Moreover, with $\pi^\ddag$ as in Lemma \ref{daggers-construction-lemma}, the permutations $\pi,\pi^\dag,\pi^\ddag$ are involutions when restricted to the set $\mathcal{I}=\{I^\dag,J^\dag,K^\dag,L^\dag,\pi(I^\dag),\pi(J^\dag),\pi(K^\dag),\pi(L^\dag)\}$,
and agree on the complement $\mathcal{I}^c$.
\end{lemma}
\begin{proof}
If we prove $\pi^\dag$ has the cycles $(I^\dag,K^\dag)$ and $(J^\dag,L^\dag)$, then from (\ref{def-pi-dag}) and $\pi^\ddag=\pi \alpha^{\pi^\dag}_{I^\dag,J^\dag}$ we have that $\pi$, $\pi^\dag$, $\pi^\ddag$ all agree on the complement of $\mathcal{I}$,
which proves the last claim in the lemma. Note that $\pi$ maps $\mathcal{I}$ onto itself and is an involution when restricted to $\mathcal{I}$. Therefore $\pi$ is an involution when restricted to $\mathcal{I}^c$, and hence so are $\pi^\dag$ and $\pi^\ddag$. So, we only need to prove $\pi^\dag$ has the cycles as claimed.

Since $\pi^\ddag=\pi \alpha^{\pi^\dag}_{I^\dag,J^\dag}$, it suffices now to show that $\pi^\dag$ is an involution on $\mathcal{I}$ and has cycles $(I^\dag,K^\dag)$ and $(J^\dag,L^\dag)$, which can be achieved by examining the cases in (\ref{def-pi-dag}) one by one. For instance, in case one, where
$\pi(I^\dagger)=K^\dagger$ and $\pi(J^\dagger) \not = L^\dagger$  we have
$\mathcal{I}=\{I^\dag,J^\dag,K^\dag,L^\dag,\pi(J^\dag),\pi(L^\dag)\}$ and
\beas
\pi^\dag(I^\dag)&=&\pi\alpha^\pi_{J^\dag,L^\dag}(I^\dag)= \pi(I^\dag)=K^\dag\\ \pi^\dag(J^\dag)&=&\pi\alpha^\pi_{J^\dag,L^\dag}(J^\dag)=\pi \tau_{J^\dag,\pi(L^\dag)}(J^\dag)= \pi(\pi(L^\dag))=L^\dag \quad \mbox{and}\\
\pi^\dag(\pi(J^\dag))&=&\pi\alpha^\pi_{J^\dag,L^\dag}\pi(J^\dag)=\pi\tau_{J^\dag,\pi(L^\dag)}\tau_{L^\dag,\pi(J^\dag)}\pi(J^\dag)=\pi(L^\dag).
\enas
Hence $\pi^\dag$ is an involution on $\mathcal{I}$ and has cycles $(I^\dag,K^\dag)$, $(J^\dag,L^\dag)$, $(\pi(J^\dag),\pi(L^\dag))$.

In case ten, $|\mathcal{I}|=8$, and $\pi^\dag$ will be an involution on $\mathcal{I}$
with cycles $(I^\dag,K^\dag)$, $(J^\dag,L^\dag)$, $(\pi(I^\dag),\pi(K^\dag))$, $(\pi(J^\dag),\pi(L^\dag))$. As an illustration we note
\beas
\pi^\dag(I^\dag)=\pi\alpha^\pi_{I^\dag,K^\dag}\alpha^\pi_{J^\dag,L^\dag}(I^\dag)=
\pi\alpha^\pi_{I^\dag,K^\dag}(I^\dag)=\pi \tau_{I^\dag,\pi(K^\dag)} \tau_{K^\dag,\pi(I^\dag)} (I^\dag)= \pi \tau_{I^\dag,\pi(K^\dag)} (I^\dag) = \pi(\pi(K^\dag))=K^\dag.
\enas
That $\pi^\dag$ has the other cycles as claimed can be shown similarly. So in these two cases $\pi^\dag$ is an involution restricted to $\mathcal{I}$ and has the cycles $(I^\dag,K^\dag)$, $(J^\dag,L^\dag)$.

That $\pi^\dag$ is an involution on $\mathcal{I}$ with cycles $(I^\dag,K^\dag)$, $(J^\dag,L^\dag)$ can be similarly shown for the other cases, completing the proof.
\end{proof}
Henceforth we will only write  $\alpha_{i,j}$ for $\alpha^\pi_{i,j}$ unless otherwise mentioned. The utility of the construction (\ref{def-pi-dag}) as a coupling is indicated by the following result.
\begin{theorem}
Suppose $\pi$ is chosen uniformly at random from $\Pi_n$ and $(I^\dag,J^\dag,K^\dag\linebreak, L^\dag)$ has joint distribution $p(\cdot)$ as in (\ref{def-p}). If $\pi^\dag$ is obtained from $\pi$ according to (\ref{def-pi-dag}) above, then $\pi$ and $\pi^\dag$ are constructed on a common space and $\pi^\dag$ satisfies the conditions specified in Lemma \ref{daggers-construction-lemma}.
 \label{thm-sqbias}
\end{theorem}
\begin{proof}
By hypothesis the indices $(I^\dagger, J^\dagger, K^\dagger, L^\dagger)$ have distribution in (\ref{def-p}). From Lemma \ref{lemma-pi-dag-involution}, we see that $\pi^\dag$ is an involution and has cycles $(I^\dag,K^\dag)$, $(J^\dag,L^\dag)$.
It only remains to verify that the distribution of $\pi^\dagger$ is uniform over all involutions in $\Pi_n$ having cycles $(I^\dagger,K^\dagger)$ and $(J^\dagger,L^\dagger)$. That is, recalling that
${\bf I}^\dagger=(I^\dagger,J^\dagger)$ and ${\bf K}^\dagger=(K^\dagger, L^\dagger)$, letting
$$
\Pi_{n,{\bf I}^\dagger, {\bf K}^\dagger}= \{ \pi \in \Pi_n: \pi(I^\dagger)=K^\dagger, \pi(J^\dagger)=L^\dagger\},
$$
 we need to verify that
\bea \label{for-phi}
P(\pi^\dag=\phi|{\bf I}^\dagger, {\bf K}^\dagger)=\frac{1}{|\Pi_{n,{\bf I}^\dagger, {\bf K}^\dagger}|}=\frac{1}{|\Pi_{n-4}|} \quad \mbox{for all $\phi \in \Pi_{n,{\bf I}^\dagger, {\bf K}^\dagger}$.}
\ena

\par
Since $I^\dagger,J^\dagger,K^\dagger, L^\dagger$ are distinct, with $\mathcal{I}$ as in Lemma \ref{lemma-pi-dag-involution}, the size of $\mathcal{I}$ satisfies
$4 \le |\mathcal{I}| \le 8$. In addition, since $\pi$ is an involution we see that $|\{I^\dag,J^\dag,K^\dag,L^\dag\}\cap\{\pi(I^\dag),\pi(J^\dag),\pi(K^\dag),\pi(L^\dag)\}|$
is even. Hence so is $|\mathcal{I}|$, and we conclude that $\mathcal{I} \in \{4,6,8\}$.
\par

For $\phi \in \Pi_{n,{\bf I}^\dagger, {\bf K}^\dagger}$, independence of $\mathcal{I}$ and $(\mathbf{I}^\dag,\mathbf{K}^\dag)$ yields
\bea
 P(\pi^\dag=\phi|{\bf I}^\dagger, {\bf K}^\dagger)&=&\sum_{\iota \in \{4, 6, 8\}} P\left(\left.\pi^\dag=\phi\right|\,|\mathcal{I}|=\iota,{\bf I}^\dagger, {\bf K}^\dagger\right)P(\left.|\mathcal{I}|=\iota\,\right|{\bf I}^\dagger, {\bf K}^\dagger)\nonumber\\
&=&\sum_{\iota \in \{4, 6, 8\}}P\left(\left.\pi^\dag=\phi\right|\,|\mathcal{I}|=\iota,{\bf I}^\dagger, {\bf K}^\dagger\right)P(|\mathcal{I}|=\iota).
\label{sum-cond-pis}
\ena

For $\pi \in \Pi_n$ let $\overline{\pi}$ denote the restriction of $\pi$ to the complement of $I^\dagger,J^\dagger,K^\dagger,L^\dagger$. First consider the case $\iota=4$ that is $\mathcal{I}=\{I^\dag,J^\dag,K^\dag,L^\dag\}$. Since $\pi^\dagger \in \Pi_{n,{\bf I}^\dagger,{\bf K}^\dagger}$ the permutation $\pi^\dagger$ agrees with every $\phi \in \Pi_{n,{\bf I}^\dagger,{\bf K}^\dagger}$ on $I^\dagger,J^\dagger,K^\dagger,L^\dagger$, and, as
$\pi^\dagger$ and $\pi$ agree on $\mathcal{I}^c$,
\beas
P\left(\left.\pi^\dag=\phi\right|\,|\mathcal{I}|=4,{\bf I}^\dagger, {\bf K}^\dagger\right) &=&P\left(\left.\overline{\pi^\dag}=\overline{\phi}\right|\,|\mathcal{I}|=4,{\bf I}^\dagger, {\bf K}^\dagger\right)\\
&=&P\left(\left.\overline{\pi}=\overline{\phi}\right|\,|\mathcal{I}|=4,{\bf I}^\dagger, {\bf K}^\dagger\right)\\
&=&\frac{1}{|\Pi_{n-4}|}.
\enas

Now suppose $\iota=6$. In this case, the set $\mathcal{J}=\mathcal{I} \setminus \{I^\dagger,J^\dagger,K^\dagger, L^\dagger\}$ has size 2, say ${\cal J}=\{i_1,i_2\}$. We claim that $\pi^\dagger$ has the cycle
$(i_1,i_2)$ with $\{i_1,i_2\}=\{\pi(a),\pi(b)\}$ for some $a,b \in \{I^\dag,J^\dag,K^\dag,L^\dag\}$, and that conditional on $|\mathcal{I}|=6$ and $\{I^\dag,J^\dag,K^\dag,L^\dag\}$ the values $i_1,i_2$ are uniform over all pairs of distinct values in $\{I^\dagger,J^\dagger,K^\dagger, L^\dagger\}^c$. That $\pi^\dag$ has the cycle $(i_1,i_2)$ follows from Lemma \ref{lemma-pi-dag-involution}. Suppose $(i_1,i_2)=(\pi(J^\dag),\pi(L^\dag))$ as in case 1 in (\ref{def-pi-dag}), since $J^\dag,L^\dag$ do not form a cycle, and their images under $\pi$ under case 1 is constrained exactly to lie outside $\{I^\dagger,J^\dagger,K^\dagger, L^\dagger\}$, as $\pi$ is uniform over $\Pi_n$, these images are uniform over $\{I^\dagger,J^\dagger,K^\dagger, L^\dagger\}^c$.
One can show that these properties hold similarly in the remaining cases.
The remaining cycles of $\pi^\dagger$ are in $\mathcal{I}^c$ and thus are the same as those of $\pi|_{\mathcal{I}^c}$. Thus, the cycles of $\overline{\pi^\dag}$ are conditionally uniform, that is, for $\phi \in \Pi_{n,{\bf I}^\dagger, {\bf K}^\dagger}$
\beas
P\left(\left.\pi^\dag=\phi\right|\,|\mathcal{I}|=6,{\bf I}^\dagger, {\bf K}^\dagger\right)
=P\left(\left.\overline{\pi^\dag}=\overline{\phi}\right|\,|\mathcal{I}|=6,{\bf I}^\dagger, {\bf K}^\dagger\right) = \frac{1}{|\Pi_{n-4}|}.
\enas

The case $\iota=8$, that is, where $R_1=R_2=0$, is handled similar to $\iota=6$. Here $\pi^\dag$ will have the cycles $(\pi(I^\dag),\pi(K^\dag)),(\pi(J^\dag),\pi(L^\dag))$. These two cycles are both of the form $(\pi(a),\pi(b))$ with $a,b\in\{I^\dag,J^\dag,K^\dag,L^\dag\}$ and hence, as in the case $\iota=6$, uniform random transpositions. Since, $\pi$ and $\pi^\dag$ agree on $\mathcal{I}^c$, we can see that $\pi^\dag$ has uniform random transpositions on $\{I^\dag,J^\dag,K^\dag,L^\dag\}^c=\mathcal{I}^c\cup \{\pi(I^\dag),\pi(J^\dag),\pi(K^\dag),\pi(L^\dag)\}$ yielding
$$
P\left(\left.\pi^\dag=\phi\right|\,|\mathcal{I}|=8,{\bf I}^\dagger, {\bf K}^\dagger\right)
=P\left(\left.\overline{\pi^\dag}=\overline{\phi}\right|\,|\mathcal{I}|=8,{\bf I}^\dagger, {\bf K}^\dagger\right) = \frac{1}{|\Pi_{n-4}|}.$$ Thus (\ref{sum-cond-pis}) now yields
\beas
P(\pi^\dag=\phi)=\frac{1}{|\Pi_{n-4}|},
\enas
verifying (\ref{for-phi}) and proving the theorem.
\end{proof}
So the tuple $(\pi^\dag,\pi^\ddag)$ obtained in Theorem \ref{thm-sqbias} satisfy the conditions in Lemma \ref{daggers-construction-lemma}. Hence $(W^\dag,W^\ddag)$ constructed from $(\pi^\dag,\pi^\ddag)$ as in Lemma \ref{daggers-construction-lemma} has the required square bias distribution.
\par We conclude this section with the calculation of the normalization constant for the distribution $p(\cdot)$ in (\ref{def-p}).
\begin{lemma}
\label{lemma-cn} For $n \ge 4$ and $D=((d_{ij}))_{n\times n}$ satisfying  (\ref{D-conditions}), we have \beq \label{def-cn}
c_n\sum_{|\{i,j,k,l\}|=4} [d_{ik}+d_{jl}-(d_{ij}+d_{kl})]^2=1 \quad
\mbox{where} \quad c_n=\frac{1}{2(n-1)^2(n-3)}=O(\frac{1}{n^3}), \eeq and in particular
\bea \label{n-ge-9}
 c_n\le\frac{1}{n^3}\quad\mbox{when
$n\ge 9$}. \ena
\end{lemma}

\begin{proof}
From (\ref{bref}), we have
 $$
 W-W'=2(d_{I\pi(I)}+d_{J\pi(J)}-(d_{IJ}+d_{\pi(I)\pi(J)})),
 $$
where $(I,J)$ are two distinct indices selected uniformly from $\{1,2,\ldots,n\}$. Since $\pi$ is an involution chosen uniformly, we have
\bea
E(W-W')^2&=&\frac{1}{n(n-1)}\sum_{i\neq j} 4E(d_{i\pi(i)}+d_{j\pi(j)}-(d_{ij}+d_{\pi(i)\pi(j)}))^2\nonumber\\
&=&\frac{4}{n(n-1)^2(n-3)}\sum_{|\{i,j,k,l\}|=4}(d_{ik}+d_{jl}-(d_{ij}+d_{kl}))^2.\label{diff2-inv}
\ena
Using (\ref{diff2-inv}), (\ref{2lambdasigma2}) and (\ref{lambda-inv}) and $\sigma_D^2=1$, we obtain
\beas
\frac{4}{n(n-1)^2(n-3)}\sum_{|\{i,j,k,l\}|=4}(d_{ik}+d_{jl}-(d_{ij}+d_{kl}))^2=2\lambda\sigma_D^2=\frac{8}{n}.\label{eq-2lambdas2}
\enas
On simplification, we obtain
\beas
\frac{1}{2(n-1)^2(n-3)}\sum_{|\{i,j,k,l\}|=4}(d_{ik}+d_{jl}-(d_{ij}+d_{kl}))^2=1,
\enas
proving (\ref{def-cn}). The verification of (\ref{n-ge-9}) is direct.
\end{proof}

\section{$L^1$ bounds}
In this section we derive the $L^1$ bounds for the normal
approximation of $W=\sum_{i=1}^n d_{i\pi(i)}$, where $\pi$ is chosen uniformly at random from $\Pi_n$. The main theorem in this section is the following.
\begin{theorem}
Let $\pi$ be an involution chosen uniformly at random from $\Pi_n$ and $D=((d_{ij}))$ be an array satisfying (\ref{D-conditions}). Then with $\beta_D$ as in (\ref{def-beta}), $W=\sum_{i=1}^n d_{i\pi(i)}$
satisfies
$$ ||F_W-\Phi||_1\le \frac{\beta_D}{n}\left( 224+1344\frac{1}{n}+384\frac{1}{n^2}\right) \quad\mbox{when  $ n\ge 9$.}$$
In particular,
$$ ||F_W-\Phi||_1\le 379\frac{\beta_D}{n}\quad\mbox{when  $ n\ge 9$.}$$
\label{L1-theorem}
\end{theorem}

\par We will need the following inequalities in order to prove Theorem \ref{L1-theorem}. To avoid writing down the indices over which we are summing, unless otherwise specified the summation will be taken over the same index set as the one immediately preceding it.
\par With $p(\cdot)$ as in (\ref{def-p}), in what follows, we will apply bounds such as
\beq
\sum_{|\{i,j,k,l\}|=4} |d_{ik}|p(i,j,k,l)&=&c_n\sum |d_{ik}|[d_{ik}+d_{jl}-(d_{ij}+d_{kl})]^2\nonumber\\
&\le & \sum_{i,j,k,l}|d_{ik}|[d_{ik}+d_{jl}-(d_{ij}+d_{kl})]^2\nonumber\\
&=& c_n\sum |d_{ik}|(d^2_{ik}+d^2_{jl}+d^2_{ij}+d^2_{kl})\nonumber\\
&\le& 4c_n n^2\beta_D\le
4\frac{\beta_D}{n}\quad\mbox{when $n\ge 9$, from
(\ref{n-ge-9})}.\label{ineq-main}\eeq The first nontrivial equality above uses the special form of the term inside squares. Whenever we encounter a cross term we always get a free index to sum over which gives us zero since $d_{i+}=0\,\forall i$. The second inequality uses the fact that for
any choices $\iota_1,\iota_2,\kappa_1,\kappa_2\in\{i,j,k,l\}$ with
$\iota_1\not=\kappa_1$ and $\iota_2\not=\kappa_2$, perhaps by
relabelling the indices after the inequality, \beas \sum_{i,j,k,l}
|d_{\iota_1\kappa_1}|d^2_{\iota_2\kappa_2}\le (\sum
|d_{ij}|^3)^\frac{1}{3}(\sum |d_{kl}|^3)^\frac{2}{3}\le
n^2\beta_D.\label{ineq-holder} \enas Generally the exponent of $n$ in
such an inequality will be 2 less the number of  indices over which
we are summing up. For instance, if we are summing up over 5 indices
the exponent of $n$ will be 3 and so on. In particular, \beq
\sum_{|\{i,j,k,l,s\}|=5} |d_{\iota s}|p(i,j,k,l)&\le& 4n^3c_n\beta_D\quad\mbox{where $\iota\in\{i,j,k,l\}$}\nonumber\\
&\le & 4\beta_D\quad\mbox{when $n\ge 9$, using (\ref{n-ge-9}) and}
\label{ineq-main2}\\
\sum_{|\{i,j,k,l,s,t\}|=6} |d_{st}|p(i,j,k,l)&\le& 4 n^4c_n\beta_D\le 4n\beta_D\quad\mbox{when $n\ge 9$}\label{ineq-main3}.
\eeq

\begin{theorem}\label{thm-w-wstar}
Suppose $D=((d_{ij}))$ is an array satisfying (\ref{D-conditions}) and $\pi$ and $\pi^\dag$ are as in Theorem \ref{thm-sqbias}, and
$\pi^\ddag, W^\dag,W^\ddag$ and $W^*$ are as in Lemma \ref{daggers-construction-lemma}. Then $W,W^\dag,W^\ddag$ can be decomposed as
\bea \label{composition}
W=S+T& W^\dag=S+T^\dag & W^\ddag=S+T^\ddag,
\ena
where
\beq
S=\sum_{i\notin\mathcal{I}} d_{i\pi(i)},\quad T=\sum_{i\in\mathcal{I}}
d_{i\pi(i)}, \quad T^\dagger= \sum_{i\in\mathcal{I}}
d_{i\pi^\dagger(i)}\quad\mbox{and}\quad T^\ddagger=\sum_{i\in\mathcal{I}}
d_{i\pi^\ddagger(i)},\label{composition-def}\eeq
where $\mathcal{I}$ is as in Lemma \ref{lemma-pi-dag-involution}.
Also, $W^*$ has the $W$ zero bias distribution and satisfies
\beq
E|W-W^*|\le 112\frac{\beta_D}{n}+672\frac{\beta_D}{n^2}+192\frac{\beta_D}{n^3}\quad\mbox{when $n\ge 9$}.
\eeq
\end{theorem}
In view of (\ref{l1result}) Theorem \ref{thm-w-wstar} implies Theorem \ref{L1-theorem}.

\begin{proof} Lemma \ref{daggers-construction-lemma} guarantees that $W^*$ has
the $W$-zero biased distribution. Recalling $\mathcal{I}=\{I^\dag,J^\dag,K^\dag,L^\dag,\linebreak \pi(I^\dag),\pi(J^\dag),\pi(K^\dag),\pi(L^\dag)\}$ and that $\pi,\pi^\dag$ and $\pi^\ddag$ agree on $\mathcal{I}^c$ by Lemma \ref{lemma-pi-dag-involution}, we obtain decomposition (\ref{composition}).

\par From $W^*=UW^\dag+(1-U)W^\ddag$ and (\ref{composition}), we obtain
\beq E|W-W^*|=E|UT^\dagger+(1-U)T^\ddagger-T|. \eeq
Using the fact that $E(U)=1/2$, and that $U$ is independent of $T^\dag$ and $T^\ddag$, we obtain
\beq
E|W-W^*|\le \frac{1}{2}(E|T^\dag|+E|T^\ddag|)+E|T|=EV \quad \mbox{where} \quad V=|T^\dag|+|T|,\label{def-v-1}
\eeq
where the equality follows from the fact that $\pi^\dag,\pi^\ddag$, and therefore $T^\dag,T^\ddag$, are exchangeable.

Thus our goal is to bound the $L^1$ norms of $T$ and $T^\dag$ and we proceed in a case by case basis, much along the lines of Section 6 in \cite{goldstein}. In summary,
we group the ten cases in (\ref{def-pi-dag}) into the following five cases: $R_1=1$; $R_1=0,R_2=1$; $R_1=2$; $R_1=0,R_2=2;R_1=0,R_2=0$.\\
\noindent\textbf{Computation on $R_1=1$:} The event $R_1=1$, which
we indicate by $\mathbf{1}_1$, can occur in four different ways, corresponding to the first four cases in the definition of $\pi^\dag$ in (\ref{def-pi-dag}). With $V$ as in (\ref{def-v-1}), we can decompose $\mathbf{1}_1$ to yield
\beq
V\mathbf{1}_1=V\mathbf{1}_{1,1}+V\mathbf{1}_{1,2}+V\mathbf{1}_{1,3}+V\mathbf{1}_{1,4},\label{decomp1}
\eeq
where ${\bf 1}_{1,1}={\bf 1}(\pi(I^\dag)=K^\dag,\pi(J^\dag)\not=L^\dag)$ and ${\bf 1}_{1,m}$ for $m=2,3,4$ similarly corresponding to the other three cases in (\ref{def-pi-dag}), in their respective order. On $\mathbf{1}_{1,1}$, we have $\mathcal{I}=\{I^\dag,J^\dag,K^\dag,L^\dag,\pi(J^\dag),\pi(L^\dag)\}$ and $\pi(I^\dag)=K^\dag$, yielding
\beq
T\mathbf{1}_{1,1}&=&\sum_{i\in\mathcal{I}} d_{i\pi(i)}\mathbf{1}_{1,1}= 2(d_{I^\dag K^\dag}+d_{J^\dag\pi(J^\dag)}+d_{L^\dag\pi(L^\dag)})\mathbf{1}_{1,1}.\label{def-t-1.1}
\eeq
By Lemma \ref{lemma-pi-dag-involution} $\pi^\dag$ has cycles $(I^\dag,K^\dag),(J^\dag,L^\dag)$ and is an involution restricted to $\mathcal{I}$,
hence
\beq
T^\dag\mathbf{1}_{1,1}&=&\sum_{i\in\mathcal{I}} d_{i\pi^\dag(i)}\mathbf{1}_{1,1}=2(d_{I^\dag K^\dag}+d_{J^\dag L^\dag}+d_{\pi(J^\dag)\pi(L^\dag)})\mathbf{1}_{1,1}. \label{def-1.1}
\eeq
So, we obtain
\beq
E|T|\mathbf{1}_{1,1}& \le & 2(E|d_{I^\dag K^\dag}|\mathbf{1}_{1,1}+E|d_{J^\dag\pi(J^\dag)}|\mathbf{1}_{1,1}+E|d_{L^\dag\pi(L^\dag)}|\mathbf{1}_{1,1})\label{bd-T-1.1}\\
E|T^\dag|\mathbf{1}_{1,1}&\le& 2(E|d_{I^\dag,K^\dag}|\mathbf{1}_{1,1}+E|d_{J^\dag,L^\dag}|\mathbf{1}_{1,1}+E|d_{\pi(J^\dag),\pi(L^\dag)}|\mathbf{1}_{1,1}).\label{qtytobound}
\eeq
Because of the indicator in (\ref{bd-T-1.1}) and (\ref{qtytobound}), we need to consider the joint distribution
\beas
p_2(i,j,k,l,s,t)=P((I^\dag,J^\dag,K^\dag,L^\dag,\pi(I^\dag),\pi(J^\dag))=(i,j,k,l,s,t)),
\enas
which includes the images of $I^\dag,J^\dag$ under $\pi$, say $s$ and $t$ respectively. With $c_n$ as in Lemma \ref{lemma-cn}, we claim $p_2(\cdot)$ is given by
\beq\label{def-p2}
p_2(i,j,k,l,s,t)=\left\{\begin{array}{clcr}
\frac{1}{n-1}p(i,j,k,l)&\mbox{when $s=j,t=i$}\\
0&\mbox{when $s=j,t\not=i$ or $s\not=j,t=i$}\\
0&\mbox{when $s=i$ or $t=j$}\\
0&\mbox{when $s=t$}\\
\frac{1}{(n-1)(n-3)}p(i,j,k,l)&
\mbox{when $t\notin\{j,s,i\} \Leftrightarrow s\not\in \{i,t,j\}$.}\\
\end{array}
\right.
\eeq
To justify (\ref{def-p2}), note first that $s=j$ if and only if $t=i$, for example, $s=j$ implies $t=\pi(j)=\pi(s)=\pi(\pi(i))=i$, and therefore
\beq
\{s=j\}=\{s=j,t=i\}=\{t=i\}.\label{reduction}
\eeq
Thus the second case of (\ref{def-p2}) has zero probability. The remaining trivial
cases can be discarded using the fact that $\pi \in \Pi_n$.
Leaving these out, the first probability is derived using (\ref{reduction}), since
the image of $I^\dagger$ under $\pi$ is uniform over $\{I^\dagger\}^c$ and independent of $(I^\dag,J^\dag,K^\dag,L^\dag)$, so in particular takes the value $s=j$ with probability $1/(n-1)$. In the last case it is easy to see that $t\notin\{j,s,i\}$ and $s\not\in \{i,t,j\}$ are
equivalent. The image of $I^\dagger$ is uniform over all available $n-1$ choices, and, conditional on $\pi(I^\dagger) \not = J^\dagger$, the $n-3$ remaining choices for the image of $J^\dag$ fall in $\{I^\dagger,\pi(I^\dagger),J^\dagger\}^c$ uniformly.

Next we bound each of the summands in (\ref{bd-T-1.1}) and (\ref{qtytobound}) separately. First note that under $\mathbf{1}_{1,1}$ only the last form of
$p_2(i,j,k,l,s,t)$ in (\ref{def-p2}) is relevant. In particular, $s \not = t$ always, and since $i,j,k,l$ are distinct, $s=k$ implies $s \not \in \{i,j\}$. Hence, for the first summand in (\ref{bd-T-1.1}), we obtain
\beq
E|d_{I^\dag K^\dag}|\mathbf{1}_{1,1}&=&\sum_{i,j,k,l,s,t} |d_{ik}|p_2(i,j,k,l,s,t)\mathbf{1}(s=k,t\not=l)\nonumber\\
&=&\sum_{|\{i,j,k,l,t\}|=5} |d_{ik}|p_2(i,j,k,l,k,t)\nonumber\\
&= &\frac{n-4}{(n-1)(n-3)}\sum_{|\{i,j,k,l\}|=4}|d_{ik}|p(i,j,k,l)\nonumber\\
&\le &\frac{4}{n(n-1)}\beta_D\quad\mbox{using (\ref{ineq-main})}\nonumber\\
&\le & 8\frac{\beta_D}{n^2}\quad\mbox{when $n\ge 9$}.\label{bd-i-k-1.1}
\eeq
Similarly, we can estimate the second summand in (\ref{bd-T-1.1}) as
\beq
E|d_{J^\dag\pi(J^\dag)}|\mathbf{1}_{1,1}&=& \sum_{i,j,k,l,s,t} |d_{jt}|p_2(i,j,k,l,s,t)\mathbf{1}(s=k,t\neq l)\nonumber\\
&=&\sum_{|\{i,j,k,l,t\}|=5} |d_{jt}|p_2(i,j,k,l,k,t)\nonumber\\
&= &\frac{1}{(n-1)(n-3)}\sum |d_{jt}|p(i,j,k,l)\nonumber\\
&\le & \frac{4}{(n-1)(n-3)}\beta_D\quad\mbox{using (\ref{ineq-main2}), when $n\ge 9$}\nonumber\\
&\le & 8\frac{\beta_D}{n^2}\quad\mbox{when $n\ge 9$}.\label{bd-j-pi-j-1.1}
\eeq
The last summand in (\ref{bd-T-1.1}) is $E|d_{L^\dag\pi(L^\dag)}|\mathbf{1}_{1,1}$. Using the fact that $$(I^\dag,J^\dag,K^\dag,L^\dag,\pi(I^\dag),\pi(J^\dag))\stackrel{\mathcal{L}}{=}(K^\dag,L^\dag,I^\dag,J^\dag,\linebreak\pi(K^\dag),\pi(L^\dag)),$$ where $\stackrel{\mathcal{L}}{=}$ denotes equality in distribution, (\ref{bd-j-pi-j-1.1}) yields,
\beq
E|d_{L^\dag\pi(L^\dag)}|\mathbf{1}_{1,1}=E|d_{J^\dag\pi(J^\dag)}|\mathbf{1}_{1,1}\le
8\frac{\beta_D}{n^2}\quad\mbox{when $n\ge 9$}.\label{bd-l-pi-l-1.1}
\eeq
Thus  combining the bounds in (\ref{bd-i-k-1.1}),(\ref{bd-j-pi-j-1.1}),(\ref{bd-l-pi-l-1.1}) we obtain
the following bound on the term (\ref{bd-T-1.1}),
\beq
E|T|\mathbf{1}_{1,1}\le 48\frac{\beta_D}{n^2}\quad\mbox{when $n\ge 9$}.\label{bd-t-1.1}
\eeq

Now moving to (\ref{qtytobound}), we note the first summand in (\ref{qtytobound}) is the same as the first summand of (\ref{bd-T-1.1}), and so can be bounded by (\ref{bd-i-k-1.1}). We can bound $E|d_{J^\dag,L^\dag}|\mathbf{1}_{1,1}$, which is the second summand in (\ref{qtytobound}) in a similar fashion as in (\ref{bd-i-k-1.1}) through the following calculation
\beq
E|d_{J^\dag L^\dag}|\mathbf{1}_{1,1}&=&\sum_{|\{i,j,k,l,t\}|=5} |d_{jl}|p_2(i,j,k,l,k,t)\nonumber\\
&=&\frac{n-4}{(n-1)(n-3)}\sum_{|\{i,j,k,l\}|=4}|d_{jl}|p(i,j,k,l)\nonumber\\
&\le & 8\frac{\beta_D}{n^2}\quad\mbox{when $n\ge 9$, using (\ref{ineq-main}).}\label{bd-j-l-1.1}
\eeq

 So we are only left with the last summand of (\ref{qtytobound}) which is $E|d_{\pi(J^\dag)\pi(L^\dag)}|\mathbf{1}_{1,1}$. For this we will need to introduce the joint distribution \bea \label{intro-p3}p_3(i,j,k,l,s,t,r)=P((I^\dag,J^\dag,K^\dag,L^\dag,\pi(I^\dag),\pi(J^\dag),\pi(L^\dag))=(i,j,k,l,s,t,r)).\ena
The case $\mathbf{1}_{1,1}$ is equivalent to $s=k,t\not=l,r\neq j$, and, since $\pi\in\Pi_n$ it is equivalent to $s=k$ and $|\{i,j,k,l,r,t\}|=6$. We claim that
\bea\label{def-p3}
p_3(i,j,k,l,s,t,r)= \frac{1}{(n-1)(n-3)(n-5)}p(i,j,k,l)&\mbox{when $s=k\,\mbox{and}\,|\{i,j,k,l,r,t\}|=6$}.
\ena
 To justify (\ref{def-p3}), we note that since $\pi$ is independent of $\{I^\dag,J^\dag,K^\dag,L^\dag\}$, the image of $I^\dag$ is uniform over $n-1$ choices from $\{I^\dag\}^c$ and conditional on $\pi(I^\dag)=K^\dag$, the $n-3$ choices for $\pi(J^\dag)$ fall in $\{I^\dag,J^\dag,K^\dag\}^c$ uniformly. Conditioned on these two images, when $\pi(J^\dag)\neq L^\dag$, $\pi(L^\dag)$ is distributed uniformly over the $n-5$ choices of $\{I^\dag,J^\dag,K^\dag,L^\dag,\pi(J^\dag)\}^c$.

 Now we can bound $E|d_{\pi(J^\dag)\pi(L^\dag)}|\mathbf{1}_{1,1}$ in the following way,
\beq
E|d_{\pi(J^\dag)\pi(L^\dag)}|\mathbf{1}_{1,1}&=& \sum_{i,j,k,l,s,t,r} |d_{tr}|p_3(i,j,k,l,s,t,r)\mathbf{1}(s=k,t\neq l,r\neq j)\nonumber\\
&=&\sum_{|\{i,j,k,l,t,r\}|=6} |d_{tr}| p_3(i,j,k,l,k,t,r)\nonumber\\
&=& \frac{1}{(n-1)(n-3)(n-5)}\sum |d_{tr}|p(i,j,k,l)\nonumber\\
&\le & \frac{4n}{(n-1)(n-3)(n-5)}\beta_D\quad\mbox{using (\ref{ineq-main3})}\nonumber\\
&\le & 16\frac{\beta_D}{n^2}\quad\mbox{when $n\ge 9$}.\label{bd-pi-j-pi-l-1.1}
\eeq
So, adding the bounds in (\ref{bd-i-k-1.1}),(\ref{bd-j-l-1.1}) and (\ref{bd-pi-j-pi-l-1.1}), and using (\ref{qtytobound}) we obtain,
\beq
E|T^\dag|\mathbf{1}_{1,1}\le 64\frac{\beta_D}{n^2}\quad\mbox{when $n\ge 9$}.\label{bd-t-dag-1.1}
\eeq
 From (\ref{bd-t-1.1}) and (\ref{bd-t-dag-1.1}), using the definition of $V$ in (\ref{def-v-1}), we obtain the following bound on the first term of (\ref{decomp1}),
 \bea
EV\mathbf{1}_{1,1}\le 112\frac{\beta_D}{n^2}\quad\mbox{when $n\ge 9$}.\label{bd-v-1.1}
\ena

Next, on $\mathbf{1}_{1,2}$, indicating the event $\pi(I^\dagger) \not = K^\dagger, \pi(J^\dagger)=L^\dagger$, we have $\mathcal{I}=\{I^\dag,J^\dag,K^\dag,L^\dag,\pi(I^\dag),\pi(K^\dag)\}$ and hence, by definition (\ref{composition-def}),
\beq
T\mathbf{1}_{1,2}&=&2(d_{I^\dag \pi(I^\dag)}+d_{K^\dag\pi(K^\dag)}+d_{J^\dag L^\dag})\mathbf{1}_{1,2}\label{def-1.2}
\eeq
We further observe, \beq
\label{1-2-exchange}
(I^\dag,J^\dag,K^\dag,L^\dag,\pi(I^\dag),\pi(J^\dag),\pi(K^\dag),\pi(L^\dag))
\stackrel{\mathcal{L}}{=}(J^\dag,I^\dag,L^\dag,K^\dag,\pi(J^\dag),\pi(I^\dag),\pi(L^\dag),\pi(K^\dag)),
\eeq
which because of (\ref{bd-t-1.1}) yields
\beq
T\mathbf{1}_{1,2}\stackrel{\mathcal{L}}{=}T\mathbf{1}_{1,1}\Rightarrow
E|T|\mathbf{1}_{1,2}=E|T|\mathbf{1}_{1,1}\le 48\frac{\beta_D}{n^2} \quad \mbox{for $n \ge 9$}\label{bd-t-1.2}
\eeq
Furthermore, the distributional equality in (\ref{1-2-exchange}) implies
\beq
T^\dag\mathbf{1}_{1,2}&=&2(d_{I^\dag K^\dag}+d_{J^\dag L^\dag}+d_{\pi(I^\dag)\pi(K^\dag)})\mathbf{1}_{1,2}\stackrel{\mathcal{L}}{=}T^\dag\mathbf{1}_{1,1},
\eeq
yielding
\beq
E|T^\dag|\mathbf{1}_{1,2}=E|T^\dag|\mathbf{1}_{1,1}\le 64\frac{\beta_D}{n^2}\quad\mbox{when $n\ge 9$}\label{bd-t-dag-1.2}.
\eeq
Thus combining (\ref{bd-t-1.2}),(\ref{bd-t-dag-1.2}) we obtain
\bea
EV\mathbf{1}_{1,2}\le 112\frac{\beta_D}{n^2}\quad\mbox{when $ n\ge 9$}.\label{bd-v-1.2}
\ena

Next on $\mathbf{1}_{1,3}$, indicating $\pi(I^\dag)=L^\dag,\pi(J^\dag)\neq K^\dag$, we have $\mathcal{I}=\{I^\dag,J^\dag,K^\dag,L^\dag,\pi(J^\dag),\pi(K^\dag)\}$ and
\beq
T\mathbf{1}_{1,3}&=&2(d_{I^\dag L^\dag}+d_{J^\dag \pi(J^\dag)}+d_{K^\dag \pi(K^\dag)})\mathbf{1}_{1,3},\label{def-t-1.3}\\
T^\dag\mathbf{1}_{1,3}&=&2(d_{I^\dag K^\dag}+d_{J^\dag L^\dag}+d_{\pi(J^\dag)\pi(K^\dag)})\mathbf{1}_{1,3}\label{def-t-dag-1.3}.
\eeq
On $\mathbf{1}_{1,3}$, we have $s=l,t\neq k$ which is equivalent to $s=l$ and $|\{i,j,k,l,t\}|=5$. Hence we may bound the first summand in (\ref{def-t-1.3}) as follows
\beq
%\label{bd-i-l-1.3}
\nonumber
E|d_{I^\dag L^\dag}|\mathbf{1}_{1,3}&=&\sum_{i,j,k,l,s,t} |d_{il}| p_2(i,j,k,l,s,t) \mathbf{1}(s=l,t\neq k)\\
&=&\sum_{|\{i,j,k,l,t\}|=5} |d_{il}|p_2(i,j,k,l,l,t)\nonumber\\
&= &\frac{n-4}{(n-1)(n-3)}\sum_{|\{i,j,k,l\}|=4} |d_{il}|p(i,j,k,l)\nonumber\\
&\le & \frac{4}{n(n-1)}\beta_D\quad\mbox{using (\ref{ineq-main})}\nonumber\\
&\le & 8\frac{\beta_D}{n^2}\quad\mbox{when $n\ge 9$.}\nonumber
\eeq
Continuing in this manner we arrive, as in (\ref{bd-v-1.2}), at
\bea
EV\mathbf{1}_{1,3}\le 112\frac{\beta_D}{n^2}\quad\mbox{when $n\ge 9$}.\label{bd-v-1.3}
\ena

Symmetries between $\mathbf{1}_{1,1}$ and $\mathbf{1}_{1,2}$ such as (\ref{1-2-exchange}) hold as well between $\mathbf{1}_{1,3}$ and $\mathbf{1}_{1,4}$, yielding
 \bea
EV\mathbf{1}_{1,4}\le 112\frac{\beta_D}{n^2}\quad\mbox{when $n\ge 9$}.\label{bd-v-1.34}
\ena

Combining the bounds from (\ref{bd-v-1.1}), (\ref{bd-v-1.2}), (\ref{bd-v-1.3}) and (\ref{bd-v-1.34}), we obtain
\beq
EV\mathbf{1}_1\le 448\frac{\beta_D}{n^2}\quad\mbox{when $n\ge 9$}.\label{bd-v-1}
\eeq

\noindent\textbf{Computation on $R_1=0,R_2=1$:} Now we need to make the following decomposition
\beas
V\mathbf{1}_2=V\mathbf{1}_{2,1}+V\mathbf{1}_{2,2},
\enas
where $\mathbf{1}_2=\mathbf{1}(R_1=0,R_2=1)$, $\mathbf{1}_{2,1}=\mathbf{1}(\pi(I^\dag)=J^\dag,\pi(K^\dag)\neq L^\dag)$, $\mathbf{1}_{2,2}=\mathbf{1}(\pi(I^\dag)\not=J^\dag,\pi(K^\dag)=L^\dag)$. On $\mathbf{1}_{2,1}$ we have $\mathcal{I}=\{I^\dag,J^\dag,K^\dag,L^\dag,\pi(K^\dag),\pi(L^\dag)\}$ which gives
\beas
T\mathbf{1}_{2,1}&=&2(d_{I^\dag J^\dag}+d_{K^\dag \pi(K^\dag)}+d_{L^\dag \pi(L^\dag)})\mathbf{1}_{2,1}\quad\mbox{and}\\
T^\dag\mathbf{1}_{2,1}&=&2(d_{I^\dag K^\dag}+d_{J^\dag L^\dag}+d_{\pi(K^\dag)\pi(L^\dag)})\mathbf{1}_{2,1}.
\enas
So, we obtain,
\beq
E|T|\mathbf{1}_{2,1}&\le & 2(E|d_{I^\dag J^\dag}|\mathbf{1}_{2,1}+E|d_{K^\dag \pi(K^\dag)}|\mathbf{1}_{2,1}+E|d_{L^\dag \pi(L^\dag)}|\mathbf{1}_{2,1}),\label{def-t-2.2}\\
E|T^\dag|\mathbf{1}_{2,1}&\le & 2(E|d_{I^\dag K^\dag}|\mathbf{1}_{2,1}+E|d_{J^\dag L^\dag}|\mathbf{1}_{2,1}+E|d_{\pi(K^\dag)\pi(L^\dag)}|\mathbf{1}_{2,1}).\label{qtytobound2}
\eeq
 Since $p(i,j,k,l)=p(i,k,j,l)$ and $\pi$ is chosen independently of $\{I^\dag,J^\dag,K^\dag,L^\dag\}$, we have
\beq (I^\dag,K^\dag,J^\dag,L^\dag,\pi(I^\dag),\pi(K^\dag),\pi(J^\dag),\pi(L^\dag))\stackrel{\mathcal{L}}{=}(I^\dag,J^\dag,K^\dag,L^\dag,\pi(I^\dag),\pi(J^\dag),\pi(K^\dag),\pi(L^\dag))\label{symmetry2}.
\eeq
Hence we obtain
$$ T\mathbf{1}_{2,1}\stackrel{\mathcal{L}}{=}T\mathbf{1}_{1,1},$$ which, by (\ref{bd-t-1.1}) gives
\bea
\label{bd-t-2.2}E|T|\mathbf{1}_{2,1}=E|T|\mathbf{1}_{1,1}\le48\frac{\beta_D}{n^2}\quad\mbox{when $n\ge 9$}.
\ena
To begin bounding $E|T^\dag|\mathbf{1}_{2,1}$, we bound the first summand in (\ref{qtytobound2}), $E|d_{I^\dag K^\dag}|\mathbf{1}_{2,1}$, as follows
\beq
E|d_{I^\dag K^\dag}|\mathbf{1}_{2,1}&=& \sum_{|\{i,j,k,l,r\}|=5} |d_{ik}|P((I^\dag,K^\dag,J^\dag,L^\dag,\pi(I^\dag),\pi(K^\dag))=(i,k,j,l,j,r))\nonumber\\
&=&\sum |d_{ik}|P((I^\dag,J^\dag,K^\dag,L^\dag,\pi(I^\dag),\pi(J^\dag))=(i,k,j,l,j,r))\quad\mbox{using (\ref{symmetry2})},\nonumber\\
&=& \sum |d_{ik}|p_2(i,k,j,l,j,r)\nonumber\\
&=& \frac{n-4}{(n-1)(n-3)}\sum_{|\{i,j,k,l\}|=4} |d_{ik}|p(i,j,k,l)\nonumber\\
&\le & 8\frac{\beta_D}{n^2}\quad\mbox{when $n\ge 9$, using (\ref{ineq-main})}.\label{bd-i-k-2.1}
\eeq
Also using the distributional equality $(I^\dag,J^\dag,K^\dag,L^\dag,\pi(I^\dag),\pi(K^\dag))\stackrel{\mathcal{L}}{=}(J^\dag,I^\dag,L^\dag,K^\dag,\pi(J^\dag),\pi(L^\dag))$ and the bound in (\ref{bd-i-k-2.1}) above, we obtain
\beq
E|d_{J^\dag L^\dag}|\mathbf{1}_{2,1}=E|d_{I^\dag K^\dag}|\mathbf{1}_{2,1}\le 8\frac{\beta_D}{n^2}\quad\mbox{when $n\ge 9$}.\label{bd-j-l-2.1}
\eeq
Using the distributional equality in (\ref{symmetry2}) and the bound in (\ref{bd-pi-j-pi-l-1.1})\beq
E|d_{\pi(K^\dag)\pi(L^\dag)}|\mathbf{1}_{2,1}=E|d_{\pi(J^\dag)\pi(L^\dag)}|\mathbf{1}_{1,1}\le 16\frac{\beta_D}{n^2}.\label{bd-pi-k-pi-l-2.1}
\eeq
Combining (\ref{bd-j-l-2.1}),(\ref{bd-pi-k-pi-l-2.1}) and using (\ref{qtytobound2}) we obtain
\beq
\label{bd-t-dag-2.2}E|T^\dag|\mathbf{1}_{2,1}\le 64\frac{\beta_D}{n^2}\quad\mbox{when $n \ge 9$}.
\eeq
Adding the two bounds in (\ref{bd-t-2.2}) and (\ref{bd-t-dag-2.2}) we obtain
\beq
EV\mathbf{1}_{2,1}\le 112\frac{\beta_D}{n^2}\quad\mbox{when $n \ge 9$}.\label{bd-2.1}
\eeq

Next, on $\mathbf{1}_{2,2}$, where $\pi(I^\dag)\not=J^\dag,\pi(K^\dag)=L^\dag$, we have $\mathcal{I}=\{I^\dag,J^\dag,K^\dag,L^\dag,\pi(I^\dag),\pi(J^\dag)\}$ and hence
\beq
T\mathbf{1}_{2,2}&=&2(d_{I^\dag\pi(I^\dag)}+d_{J^\dag\pi(J^\dag)}+d_{K^\dag L^\dag})\mathbf{1}_{2,2}\label{def-t-2.1}\\
T^\dag\mathbf{1}_{2,2}&=& 2(d_{I^\dag K^\dag}+d_{J^\dag L^\dag}+d_{\pi(I^\dag)\pi(J^\dag)})\mathbf{1}_{2,2}\label{def-t-dag-2.1}.
\eeq
Noting the distributional equality
$$
(I^\dag,J^\dag,K^\dag,L^\dag,\pi(I^\dag),\pi(J^\dag),\pi(K^\dag),\pi(L^\dag))
\stackrel{\mathcal{L}}{=}(K^\dag,L^\dag,I^\dag,J^\dag,\pi(K^\dag),\pi(L^\dag),\pi(I^\dag),\pi(J^\dag))
$$
we obtain
\beas
T\mathbf{1}_{2,1}\stackrel{\mathcal{L}}{=}T\mathbf{1}_{2,2}\quad\mbox{and}\quad T^\dag\mathbf{1}_{2,1}\stackrel{\mathcal{L}}{=}T^\dag\mathbf{1}_{2,2},
\enas
which yields
\beas
E|T|\mathbf{1}_{2,2}=E|T|\mathbf{1}_{2,1}\le 48\frac{\beta_D}{n^2}\quad\mbox{and}\quad E|T^\dag|\mathbf{1}_{2,2}=E|T^\dag|\mathbf{1}_{2,1}\le 64\frac{\beta_D}{n^2}.
\enas
Hence we have
\beq
EV\mathbf{1}_{2,2}\le 112\frac{\beta_D}{n^2}\quad\mbox{when $n\ge 9$}.\label{bd-2.2}
\eeq
Combining (\ref{bd-2.1}) with (\ref{bd-2.2}) we obtain
\beq
EV\mathbf{1}_2\le 224\frac{\beta_D}{n^2}\quad\mbox{when $n \ge 9$}.\label{bd-v-2}
\eeq

\noindent\textbf{Computation on $R_1=2$:} Here we need the decomposition
\beas
V\mathbf{1}_3=V\mathbf{1}_{3,1}+V\mathbf{1}_{3,2},
\enas
where $\mathbf{1}_3=\mathbf{1}(R_1=2)$, $\mathbf{1}_{3,1}=\mathbf{1}(\pi(I^\dag)=K^\dag,\pi(J^\dag)=L^\dag)$ and $\mathbf{1}_{3,2}=\mathbf{1}(\pi(I^\dag)=L^\dag,\pi(J^\dag)=K^\dag)$.
Note that $\mathbf{1}_{3,1}$ and $\mathbf{1}_{3,2}$ correspond to the
cases in (\ref{def-pi-dag}) where $(R_1,R_2)=(2,0)$ and $(R_1,R_2)=(2,2)$, respectively.
On both $\mathbf{1}_{3,1}$ and $\mathbf{1}_{3,2}$, we have $\mathcal{I}=\{I^\dag,J^\dag,K^\dag,L^\dag\}$ and
\beq
T^\dag\mathbf{1}_{3,1}=T\mathbf{1}_{3,1}=2(d_{I^\dag K^\dag}+d_{J^\dag L^\dag})\mathbf{1}_{3,1}\quad\mbox{since $\pi\mathbf{1}_{3,1}=\pi^\dag\mathbf{1}_{3,1}$.}\label{eq-t-tdag-3.1}
\eeq
From $(I^\dag,J^\dag,K^\dag,L^\dag,\pi(I^\dag),\pi(J^\dag))\stackrel{\mathcal{L}}{=}(J^\dag,I^\dag,L^\dag,K^\dag,\pi(J^\dag),\pi(I^\dag))$, it is clear that $E|d_{I^\dag K^\dag}|\mathbf{1}_{3,1}=E|d_{J^\dag L^\dag}|\mathbf{1}_{3,1}$ and hence it is enough to bound any one of the two summands in (\ref{eq-t-tdag-3.1}). We bound $E|d_{I^\dag K^\dag}|\mathbf{1}_{3,1}$ using (\ref{ineq-main}) as follows
\beq
E|d_{I^\dag K^\dag}|\mathbf{1}_{3,1}&=& \sum_{i,j,k,l,s,t} |d_{ik}| p_2(i,j,k,l,s,t) \mathbf{1}(s=k,t=l)\nonumber\\
&=&\sum_{|\{i,j,k,l\}|=4} |d_{ik}|p_2(i,j,k,l,k,l)\nonumber\\
&=& \frac{1}{(n-1)(n-3)}\sum|d_{ik}| p(i,j,k,l)\nonumber\\
&\le & 8\frac{\beta_D}{n^3}\quad\mbox{when $n\ge 9$}.\label{bd-i-k-3.1}\eeq
Thus, using (\ref{eq-t-tdag-3.1}) and (\ref{bd-i-k-3.1}),we obtain
\beq
EV\mathbf{1}_{3,1}\le 64\frac{\beta_D}{n^3}\quad\mbox{when $n\ge 9$}.\label{bd-3.1}
\eeq
On $\mathbf{1}_{3,2}$, we have
\bea
T\mathbf{1}_{3,2}=2(d_{I^\dag L^\dag}+d_{J^\dag K^\dag})\quad\mbox{and}\quad T^\dag\mathbf{1}_{3,2}=2(d_{I^\dag K^\dag}+d_{J^\dag L^\dag}).\label{def-case-3.2}
\ena
To obtain bounds for $E(V\mathbf{1}_{3,2})$, we bound the first summand in (\ref{def-case-3.2}) as follows,
\beq
E|d_{I^\dag L^\dag}|\mathbf{1}_{3,2}&=&\sum_{|\{i,j,k,l\}|=4} |d_{il}|p_2(i,j,k,l,l,k)\nonumber\\
&=& \frac{1}{(n-1)(n-3)}\sum |d_{il}|p(i,j,k,l)\le 8\frac{\beta_D}{n^3}\quad\mbox{when $n\ge 9$,using (\ref{ineq-main})}.\label{bd-i-l-3.2}
\eeq
Similarly its easy to obtain bounds on the other summands also and conclude as in (\ref{bd-3.1}),
\beq
EV\mathbf{1}_{3,2}\le64\frac{\beta_D}{n^3}\quad\mbox{for $n\ge 9$}.\label{bd-3.2}
\eeq
Combining (\ref{bd-3.1}) and (\ref{bd-3.2}), we obtain	
\beq
EV\mathbf{1}_3\le 128\frac{\beta_D}{n^3}\quad\mbox{when $ n\ge 9$}.\label{bd-v-3}
\eeq
\noindent\textbf{Computation for $R_1=0,R_2=2$:}
This event is indicated by $\mathbf{1}_4=\mathbf{1}(\pi(I^\dag)=J^\dag,\pi(K^\dag)=L^\dag)$. On $\mathbf{1}_4$, we have
\beq
T\mathbf{1}_4=2(d_{I^\dag J^\dag}+d_{K^\dag L^\dag})\mathbf{1}_4\label{def-t-4},\\
T^\dag\mathbf{1}_4=2(d_{I^\dag K^\dag}+d_{J^\dag L^\dag})\mathbf{1}_4\label{def-t-dag-4}.
\eeq
For the first term in (\ref{def-t-4}),
\beq
E|d_{I^\dag J^\dag}|\mathbf{1}_4 &=& \sum_{|\{i,j,k,l\}|=4} |d_{ij}| P((I^\dag,K^\dag,J^\dag,L^\dag,\pi(I^\dag),\pi(K^\dag)=(i,k,j,l,j,l))\nonumber\\
&=&\sum |d_{ij}|p_2(i,k,j,l,j,l)\quad\mbox{using (\ref{symmetry2})}\nonumber\\
& \le & \frac{1}{(n-1)(n-3)}\sum |d_{ij}|p(i,j,k,l)\nonumber\\
&\le & 8\frac{\beta_D}{n^3}\quad\mbox{when $n\ge 9$}\label{bd-i-j-4}.
\eeq
For the other summand in (\ref{def-t-4}) and (\ref{def-t-dag-4}), we can follow similar calculations as in (\ref{bd-i-j-4}) above and obtain the same bounds. Thus we will finally obtain
\beq
E|T\mathbf{1}_4|,E|T^\dag\mathbf{1}_4|\le 32\frac{\beta_D}{n^3}\quad\mbox{when $ n\ge 9$}\label{bd-t-4}.
\eeq
So we have
\beq
EV\mathbf{1}_4\le 64\frac{\beta_D}{n^3}\quad\mbox{when $n\ge 9$}.\label{bd-v-4}
\eeq
\noindent\textbf{Computation on $R_1=R_2=0$:} Now we bound the $L^1$ contribution from the last case in  (\ref{def-pi-dag}) denoted by $\mathbf{1}_5=\mathbf{1}(R_1=R_2=0)$. Here we have $\mathcal{I}=\{I^\dag,J^\dag,K^\dag,L^\dag,\pi(I^\dag),\pi(J^\dag),\pi(K^\dag),\pi(L^\dag)\}$ and
\beq
T\mathbf{1}_5=2\left(d_{I^\dag\pi(I^\dag)}+d_{J^\dag\pi(J^\dag)}+d_{K^\dag\pi(K^\dag)}+d_{L^\dag\pi(L^\dag)}\right)\label{def-t-5}\\
T^\dag\mathbf{1}_5=2\left(d_{I^\dag K^\dag}+d_{J^\dag L^\dag}+d_{\pi(I^\dag)\pi(K^\dag)}+d_{\pi(J^\dag)\pi(L^\dag)}\right)\label{def-t-dag-5}.
\eeq
Since we are on $\mathbf{1}_5$, we need to consider $p_3(\cdot)$ as introduced in (\ref{intro-p3}). On $\mathbf{1}_5$, $p_3(\cdot)$ is given by
\bea\label{def-p3-2}
p_3(i,j,k,l,s,t,r)=\frac{1}{(n-1)(n-3)(n-5)}p(i,j,k,l)\quad\mbox{when $|\{i,j,k,l,s,t,r\}|=7$}.\ena
The justification for (\ref{def-p3-2}) is essentially the same as that for (\ref{def-p3}).
\par Using $p_3(\cdot)$, we begin to bound the summands in (\ref{def-t-5}).
\beq
E|d_{I^\dag \pi(I^\dag)}|\mathbf{1}_5&=&\sum_{|\{i,j,k,l,s,t,r\}|=7} |d_{is}|p_3(i,j,k,l,s,t,r)\nonumber\\
&=& \frac{(n-6)(n-5)}{(n-1)(n-3)(n-5)}\sum_{|\{i,j,k,l,s\}|=5} |d_{is}|p(i,j,k,l)\nonumber\\
&\le & 8\frac{\beta_D}{n}\quad\mbox{when $n\ge 9$}.\label{bd-i-pi-i-5}
\eeq
It is easy to see that $I^\dag,J^\dag,K^\dag,L^\dag$ have identical marginal distributions and since $\pi$ is chosen independently of these indices, we have $E|d_{N \pi(N)}|\mathbf{1}_5$ is constant over $\{I^\dag,J^\dag,K^\dag,L^\dag\}$. Thus we obtain
\beq
E|T|\mathbf{1}_5\le 64\frac{\beta_D}{n}\quad\mbox{when $ n\ge 9$}.\label{bd-t-5}
\eeq

Now bounding the $L^1$ norm of the first summand in (\ref{def-t-dag-5}),
\beq
E|d_{I^\dag K^\dag}|\mathbf{1}_5&=&\sum_{|\{i,j,k,l,s,t,r\}|=7} |d_{ik}|p_3(i,j,k,l,s,t,r)\nonumber\\
&= & \frac{(n-4)(n-5)(n-6)}{(n-1)(n-3)(n-5)}\sum_{|\{i,j,k,l\}|=4}|d_{ik}|p(i,j,k,l)\nonumber\\
&\le & 4\frac{\beta_D}{n}\quad\mbox{when $n\ge 9$, using  (\ref{ineq-main})}.\label{bd-i-k-5}
\eeq

Now consider the last summand of (\ref{def-t-dag-5}),
\beq
E|d_{\pi(J^\dag)\pi(L^\dag)}|\mathbf{1}_5&=& \sum_{|\{i,j,k,l,s,t,r\}|=7} |d_{tr}|p_3(i,j,k,l,s,t,r)\nonumber\\
&\le&\frac{n-6}{(n-1)(n-3)(n-5)}\sum_{|\{i,j,k,l,t,r\}|=6}|d_{tr}|p(i,j,k,l)\nonumber\\
&\le & 8\frac{\beta_D}{n}\quad\mbox{when $n\ge 9$, using (\ref{ineq-main3})}.\label{bd-pi-i-pi-k-5}
\eeq

Since $(I^\dag,J^\dag,K^\dag,L^\dag,\pi(I^\dag),\pi(J^\dag),\pi(K^\dag),\pi(L^\dag))\stackrel{\mathcal{L}}{=}(J^\dag,I^\dag,L^\dag,K^\dag,\pi(J^\dag),\pi(I^\dag),\pi(L^\dag),\pi(K^\dag))$, we see that
\beq
E|d_{I^\dag K^\dag}|\mathbf{1}_5=E|d_{J^\dag L^\dag}|\mathbf{1}_5\quad\mbox{and}\quad E|d_{\pi(I^\dag)\pi(K^\dag)}|\mathbf{1}_5=E|d_{\pi(J^\dag)\pi(L^\dag)}|\mathbf{1}_5.\label{bd-eq-5}
\eeq

So, using (\ref{bd-i-k-5}),(\ref{bd-pi-i-pi-k-5}) along with (\ref{bd-eq-5}), we obtain
\beq
E|T^\dag|\mathbf{1}_5\le 48\frac{\beta_D}{n}\quad\mbox{when $n\ge 9$}.\label{bd-t-dag-5}
\eeq

Combining (\ref{bd-t-5}) and (\ref{bd-t-dag-5}), we obtain
\beq
EV\mathbf{1}_5\le 112\frac{\beta_D}{n}\quad\mbox{when $ n\ge 9$}.\label{bd-v-5}
\eeq
Combining the bounds from (\ref{bd-v-1}),(\ref{bd-v-2}),(\ref{bd-v-3}),(\ref{bd-v-4}) and (\ref{bd-v-5}), we obtain
\beq
E|W-W^*|\le EV\le 112\frac{\beta_D}{n}+672\frac{\beta_D}{n^2}+192\frac{\beta_D}{n^3}\quad\mbox{when $n\ge 9$}.\label{bd-v}
\eeq
This completes the proof of Theorem \ref{thm-w-wstar}.
\end{proof}
\section{$L^\infty$ bounds}
In this section
we will use Theorem \ref{thm-w-wstar}  obtained in the previous section to
obtain $L^\infty$ bounds using arguments similar to those in
\cite{bolt}. It is worth noting that we can use $L^1$ along with
$L^\infty$ bounds to obtain $L^p$ bounds for any $p\ge 1$ using (\ref{ineq-p}).
The main theorem of this section is the following
\begin{theorem}
\label{prop-main} Suppose we have an $n\times n$ array $D=((d_{ij}))$ satisfying $d_{ij}=d_{ji}, d_{ii}=d_{i+}=0\,\,\forall\,\,i,j$ and $\sigma^2_D=1$. If $W=\sum d_{i\pi(i)}$ where $\pi$ is an involution chosen  uniformly at random from $\Pi_n$, then for $n> 9$
  \beq \label{prop-mean-ineq}||F_W-\Phi||_\infty\le
K\frac{\beta_D}{n} \eeq Here $K=61,702,446$ is a universal constant.
\end{theorem}

Theorem \ref{prop-main} readily implies Theorem \ref{Lp-theorem} by (\ref{ineq-p}) and Theorem \ref{L1-theorem}.

\par We claim that to prove Theorem \ref{prop-main} it is enough to consider arrays with $\beta_D/n\le\epsilon_0=1/90$, and $n\ge n_0=1000$. To prove the claim, first note that
\bea\label{bd-betad}
\beta_D^\frac{1}{3}=(\sum_{1\le i,j\le n}|d_{ij}|^3)^\frac{1}{3}\ge n^{-1/3} (\sum |d_{ij}|^2)^\frac{1}{2}=n^{-1/3}\left(\frac{(n-1)(n-3)}{2(n-2)}\right)^\frac{1}{2},
\ena
which is greater than $1/2$ for $n\ge 4$.
In (\ref{bd-betad}), the first inequality follows from H\"{o}lder's inequality and the next equality follows from the fact that $\sigma^2_D=1$. So if $n<n_0$ then $\beta_D/n>1/(8n_0)$. Since $K>\max\{1/\epsilon_0,8n_0\}$ inequality (\ref{prop-mean-ineq}) holds
if $\beta_D/n>\epsilon_0$ or $n < n_0$.
\par One useful inequality that will be used repeatedly in the proof is the following
\beq
\left(|\sum_{i=1}^k a_i|\right)^3\le \left(\sum |a_i|\right)^3\le k^2\sum |a_i|^3.\label{jensen}
\eeq
\par The proof of Theorem \ref{prop-main} proceeds by first proving several auxiliary lemmas. The first lemma helps bound the error created when truncating the array as in \cite{bolt}, page 382. For the array $D=((d_{ij}))_{n\times n}$, we define
\bea
d_{ij}^\prime=d_{ij}1_{|d_{ij}|\le\frac{1}{2}}.\label{def-prime}
\ena
Letting $\Gamma=\{(i,j):|d_{ij}|>\frac{1}{2}\}$ and $\Gamma_i=\{j:(i,j)\in\Gamma\}$ we have \beq |\Gamma_i|&\le & 8\sum_j |d_{ij}|^3\label{bd-gamma-i}\quad\mbox{ and hence}\\
|\Gamma|&\le& 8\beta_D\label{bd-gamma}.
\eeq Inequality
(\ref{bd-gamma}) has the following useful consequence, that is,
 \beq
P(Y_{D^\prime}\not=Y_D)&\le& P\left(\sum_i 1_\Gamma(i,\pi(i))\ge 1\right)\nonumber\\
&\le& E(\sum_i 1_\Gamma(i,\pi(i)))\nonumber\\
&=& \frac{|\Gamma|}{n-1}<\frac{2|\Gamma|}{n}\le \frac{16\beta_D}{n}.\label{bd-d-dp}
\eeq
\begin{lemma}
\label{lemma-bd-prime}Suppose $((d_{ij}))_{n\times n}$ satisfies the conditions in Theorem \ref{prop-main}. Then with $((d^\prime_{ij}))_{n\times n}$ defined in (\ref{def-prime}), ${\beta_D} \le \epsilon_0 n$ and $n \ge n_0$, we have,
\bea \label{lemma-bd-prime-1}
|d'_{++}| \le 4\beta_D\quad\mbox{and}\quad |d'_{i+}| \le 4\sum_j|d_{ij}|^3 \le 4\beta_D,
\ena
and therefore
\beas
|\mu_{D^\prime}|&\le& 8\frac{\beta_D}{n},\\
|\sigma^2_{D^\prime}-1|&\le& 10\frac{\beta_D}{n}\quad\mbox{and}\quad\beta_{D^\prime}\le 22\beta_D.
\enas
\end{lemma}
\begin{proof} Using $d_{++}=0$ and (\ref{bd-gamma}), we have
\beq
|d'_{++}|=|\sum_{(i,j):|d_{ij}|\le 1/2} d_{ij}| =|\sum_{(i,j)\in\Gamma} d_{ij}|
\le \sum_{(i,j)\in\Gamma} |d_{ij}|\le |\Gamma|^{\frac{2}{3}}{\beta_D}^{\frac{1}{3}}\le 4\beta_D\label{ineq-5.1-1}.
\eeq
Similarly, as $d_{i+}=0$ for all $i \in \{1,\ldots,n\}$, we obtain using (\ref{bd-gamma-i}),
\bea \label{lemma-bd-prime-2}
|d^\prime_{i+}|=|\sum_{j \notin \Gamma_i} d_{ij}|=|\sum_{j \in \Gamma_i} d_{ij}|  \le |\Gamma_i|^{2/3}\left( \sum_{j \in \Gamma_i} |d_{ij}|^3\right)^{1/3} \le
4\sum_j|d_{ij}|^3\le 4\beta_D.
\ena

Now, using (\ref{mean-var-Ye}) and (\ref{ineq-5.1-1}),
\beas
|\mu_{D^\prime}|&=&\left|\frac{d^{\prime}_{++}}{n-1}\right|
\le\frac{4\beta_D}{n-1}\le 8\frac{\beta_D}{n}.
\enas

To prove the last assertion, first note that by (\ref{lemma-bd-prime-2}) we have
\beq
\sum_{i=1}^n |d^\prime_{i+}|^2\le 4\beta_D \sum_{i=1}^n |d'_{i+}|\le 16\beta_D\sum_i\sum_j|d_{ij}|^3=16\beta_D^2.\label{ineq-5.1-2}
\eeq
Similarly one can obtain
\beq\label{ineq-5.1-3}
\sum_{i=1}^n |d^\prime_{i+}|^3\le 64\beta_D^3.
\eeq
From (\ref{mean-var-Ye}) and the fact that $\sigma^2_D=1$ we obtain the following,
\beq
\left|\sigma^2_{D^\prime}-1\right|&=&\left|\frac{2}{(n-1)(n-3)}\left((n-2)\sum_{1\le i,j\le n}(d'^2_{ij}-d^2_{ij})+\frac{1}{n-1}d^{\prime2}_{++}-2\sum_{i=1}^n d^{\prime2}_{i+}\right)\right|\nonumber\\
&\le& \frac{2}{(n-1)(n-3)}\left((n-2) \sum_{(i,j)\in\Gamma}  d^2_{ij}+\frac{1}{n-1} d^{\prime2}_{++}+2\sum_{i=1}^n d^{\prime2}_{i+}\right)\nonumber \\
&\le & \frac{2}{(n-1)(n-3)}\left( (n-2)\sum_{(i,j)\in\Gamma} d^2_{ij}+16\frac{\beta_D^2}{(n-1)}+32\beta_D^2\right)\quad\mbox{using (\ref{ineq-5.1-1}) and (\ref{ineq-5.1-2})}\\
&\le& 8\frac{\beta_D}{n-1}+32\frac{\beta_D^2}{(n-1)^2(n-3)}+64\frac{\beta_D^2}{(n-1)(n-3)}\quad\mbox{since $\sum_{(i,j)\in\Gamma} d^2_{ij}\le \beta_D^\frac{2}{3}|\Gamma|^\frac{1}{3}\le 2\beta_D$}\nonumber\\
&\le& 10\frac{\beta_D}{n}<\frac{1}{9}\quad\mbox{since ${\beta_D}/{n}\le 1/90$ and $n\ge 1000$.}
\label{bd-sigma-prime}
\eeq
Hence, we obtain $|\sigma^2_{D^\prime}-1|<1/9$. Thus using (\ref{def-ghat}) with notation as in (\ref{def-beta}), and inequality (\ref{jensen}), we obtain
\beas
\beta_{D^\prime}&=&\frac{1}{\sigma_{D'}^3}\sum_{i\neq j}\left|d'_{ij}-\frac{d'_{i+}}{n-2}-\frac{d'_{+j}}{n-2}+\frac{d'_{++}}{(n-1)(n-2)}\right|^3\\
&\le& 16\frac{1}{\sigma_{D^\prime}^3}\sum_{i\not=j}\left(|d^\prime_{ij}|^3+\frac{|d^\prime_{i+}|^3}{(n-2)^3}+\frac{|d^\prime_{+j}|^3}{(n-2)^3}+\frac{|d^\prime_{++}|^3}{((n-1)(n-2))^3}\right)\\
&\le& 16\frac{1}{\sigma_{D^\prime}^3}\left(\beta_D+128\frac{\beta_D^3}{(n-2)^3}+64n^2\frac{\beta_D^3}{((n-1)(n-2))^3}\right)\\
&\le & 22\beta_D\quad\mbox{since $|\sigma^2_{D^\prime}-1|<1/9,\beta_D/n\le1/90$},
\enas
where in the second inequality we use (\ref{ineq-5.1-1}) and (\ref{ineq-5.1-3}).
\end{proof}
\begin{lemma}
\label{dtilde1}
 \par Let $D=((d_{ij}))$ be an array as in Theorem \ref{prop-main} and $D'=((d_{ij}'))$ be as in (\ref{def-prime}). Let us denote
 $$
 \widetilde{d}_{ij}=\frac{\widehat{d'_{ij}}}{\sigma_{D^\prime}}.$$
 If $\beta_D/n \le \epsilon_0$ and $n\ge n_0$, we have $|\widetilde{d}_{ij}|\le 1$.
\end{lemma}
\begin{proof}
 Because $\sigma^2_{D^\prime}>1/2$ from Lemma \ref{lemma-bd-prime}, $\sigma_{D^\prime}>2/3$. Using this inequality, definition (\ref{def-ghat}) and (\ref{lemma-bd-prime-1}) we obtain
 \beas
 |\widetilde{d}_{ij}|&\le& \frac{1}{\sigma_{D^\prime}}\left(|d^\prime_{ij}|+\frac{1}{n-2}(|d^\prime_{i+}|+|d^\prime_{+j}|)+\frac{|d^\prime_{++}|}{(n-1)(n-2)}\right)\\
&\le& \frac{3}{4}+\frac{8}{\sigma_{D^\prime}}\frac{\beta_D}{n-2}+\frac{4}{\sigma_{D^\prime}}\frac{\beta_D}{(n-1)(n-2)}\\
&<& \frac{3}{4}+16\frac{\beta_D}{n}+4\frac{\beta_D}{n}\quad\mbox{since $n\ge 1000$}\\
&=&\frac{3}{4}+20\frac{\beta_D}{n}<1\quad\mbox{since $\beta_D/n\le 1/90$}.
\enas
This completes the proof.
\end{proof}
\par For $E=((e_{ij}))_{n\times n}$ let $F_E$ be the distribution function
of $Y_E$ and
\beq \delta_E=||F_E-\Phi||_\infty.\label{def-delta-a} \eeq
\par For $\gamma>0$ let $M_n(\gamma)$ be the set of $n\times n$ matrices $E=((e_{ij}))$ satisfying
\beas
\sigma^2_E=1,e_{ij}=e_{ji},e_{ii}=e_{i+}=0\quad\forall i,j\quad\mbox{and}\quad\beta_E\le\gamma.
\enas
Let us define
\beq
\delta(\gamma,n)=\sup_{E\in M_n(\gamma)} \delta_E.\label{def-delta}
\eeq
Also, we define \beq \label{def-delta-0}
\delta^1(\gamma,n)=\sup_{E\in M_n^1(\gamma)}\delta_E\quad\mbox{where}\quad M_n^1(\gamma)=\{E\in M_n(\gamma):\sup_{i,j} |e_{ij}|\le 1\}.\eeq
\begin{lemma}\label{bd-del-del1}
When $n\ge n_0$, with $\delta_D$ and $\delta^1(\gamma,n)$ defined in (\ref{def-delta-a}) and (\ref{def-delta-0}), for all $\gamma>0$
$$
\sup\{\delta_D:D\in M_n(\gamma),\beta_D/n\le \epsilon_0\}\le \delta^1(22\gamma,n)+36\frac{\gamma}{n}.
$$
\end{lemma}
\begin{proof}
 Let $D=((d_{ij}))\in M_n(\gamma)$ with $\beta_D/n\le 1/90$. Hence, by (\ref{bd-d-dp}) and Lemmas \ref{dtilde1} and \ref{lemma-bd-prime},
\beq
\delta_D=||F_D-\Phi||_\infty&\le & \sup_t |P(Y_{D^\prime}\le t)-\Phi(t)|+P(Y_D\neq Y_{D'})\nonumber\\
&\le &\sup_t |P(Y_{D^\prime}\le t)-\Phi(t)|+\frac{16\beta_D}{n}\nonumber\\
&\le& \delta^1(22\beta_D,n)+\sup_t |\Phi(\frac{t-\mu_{D^\prime}}{\sigma_{D^\prime}})-\Phi(t)|+\frac{16\beta_D}{n}\label{bd-phi-center}.\eeq
Since $\delta^1(\gamma,n)$ is monotone in $\gamma$,
$$
\delta^1(22\beta_D,n) \le \delta^1(22\gamma,n) \quad \mbox{when $\beta_D \le \gamma$};
$$
similarly, we bound the last term as $16\beta_D/n \le 16\gamma/n$. So, we are only left with verifying that the second term is bounded by $20\gamma/n$.

From Lemma \ref{lemma-bd-prime}, we obtain $|\sigma^2_{D'}-1| \le 1/9$ yielding in particular $\sigma_{D'}\in[2/3,4/3]$. First consider
the case where $|t|\ge 8\beta_D/n$, and for a given $t$ let $t_1=(t-\mu_{D'})/\sigma_{D'}$. Since $|\mu_{D'}|\le8\beta_D/n$
by Lemma \ref{lemma-bd-prime}, $t$ and $t_1$ will be on the same side of the origin. Next, it is easy to show that for $a>0$ we have $|t\exp(-at^2/2)|\le 1/\sqrt{a}$. Hence
\bea
\left|t\exp\left(-\frac{9}{32}(t-\mu_{D'})^2\right)\right|&\le & \left|(t-\mu_{D'})\exp\left(-\frac{9}{32}(t-\mu_{D'})^2\right)\right|+|\mu_{D'}|\nonumber\\
&\le & \frac{4}{3}+|\mu_{D'}|\nonumber\\
&\le & \frac{4}{3}(1+|\mu_{D'}|).\label{ineq-t-et2}
\ena

Since $\sigma_{D'} \ge 2/3$ and Lemma \ref{lemma-bd-prime} gives $|\sigma_{D'}^2-1| \le 10 \beta_D/n$, we find that
\beas
|\sigma_{D'}-1| = \frac{|\sigma_{D'}^2-1|}{\sigma_{D'}+1} \le 10\frac{\beta_D}{n}.
\enas

Now, by the mean value theorem , $\sigma_{D'} \in [2/3,4/3]$, and Lemma \ref{lemma-bd-prime}
\beq
\left|\Phi(t_1)-\Phi(t)\right|&\le &\max\left(\phi(\frac{t-\mu_{D^\prime}}{\sigma_{D^\prime}}),\phi(t)\right)\left|\frac{t-\mu_{D^\prime}}{\sigma_{D^\prime}}-t\right|\quad\mbox{where $\phi=\Phi'$}\nonumber\\
&\le& \frac{1}{\sqrt{2\pi}} \max\left\{\exp(-\frac{9}{32}(t-\mu_{D'})^2),\exp(-\frac{t^2}{2})\right\}\left|
\frac{t(1-\sigma_{D^\prime})}{\sigma_{D^\prime}}\right|
+\frac{1}{\sqrt{2\pi}}\left|\frac{\mu_{D^\prime}}{\sigma_{D^\prime}}\right|\nonumber\\
\nonumber &\le& \frac{3}{2\sqrt{2\pi}}|\sigma_{D'}-1|\max\left\{\left|t\exp(-\frac{9}{32}(t-\mu_{D'})^2)\right|,\left|t\exp(-\frac{t^2}{2})\right|\right\}+\frac{1}{\sqrt{2\pi}}\left|\frac{\mu_{D^\prime}}{\sigma_{D^\prime}}\right|\\
&\le &\frac{3}{2\sqrt{2\pi}}|\sigma_{D'}-1|(\frac{4}{3}(1+|\mu_{D'}|))+\frac{3}{4}|\mu_{D'}|\quad\mbox{using (\ref{ineq-t-et2})}\nonumber\\
&\le & \frac{2}{\sqrt{2\pi}}|\sigma_{D'}-1|(1+|\mu_{D'}|)+\frac{3}{4}|\mu_{D'}|\label{bd-t-et2}\\
&\le & \frac{2}{\sqrt{2\pi}}10\frac{\beta_D}{n}(1+8\frac{\beta_D}{n})+6\frac{\beta_D}{n}\nonumber\\
&\le & 17\frac{\beta_D}{n}\quad\mbox{since $\beta_D/n\le1/90$} \nonumber.
\eeq

When $|t|<8\beta_D/n$, the bound is easier. Since $t_1$ lies in the interval with boundary points $3(t-\mu_{D'})/2$ and $3(t-\mu_{D'})/4$,
we have
\beq
|t_1|\le \frac{3(|t|+|\mu_{D'}|)}{2}\label{bd-t1}.
\eeq
Now, using $|t|<8\beta_D/n$, $|\mu_{D'}|<8\beta_D/n$ and (\ref{bd-t1}), we obtain
\beas
|\Phi(t_1)-\Phi(t)|&\le&\frac{1}{\sqrt{2\pi}}|t_1-t|\\
&\le &\frac{1}{\sqrt{2\pi}}(3|t|+2|\mu_{D'}|)\\
&\le & \frac{1}{\sqrt{2\pi}}40\frac{\beta_D}{n}<20\frac{\beta_D}{n}.
\enas
Thus we obtain $$\sup_t|\Phi(\frac{t-\mu_{D'}}{\sigma_{D'}})-\Phi(t)|\le 20\frac{\beta_D}{n}\le 20\frac{\gamma}{n},$$ completing the proof.
\end{proof}

In view of Lemma \ref{bd-del-del1}, to prove
Theorem \ref{prop-main} it remains only to show $\delta^1(\gamma,n)\le
{c\gamma}/{n}$ for an explicit $c$ which we eventually determine.
Hence in the following calculations we consider only arrays $D=((d_{ij}))_{n\times n}$ in $M_n^1(\gamma)$, and  so $\sup_{i,j}|d_{ij}|\le 1$.
\par We will need the following technical lemma.
\begin{lemma}\label{lemma-smaller-b}
If $D\in M_n^1(\gamma)$ and $B=((b_{ij}))_{n-l\times n-l}$ is the array formed by removing from $D$ the $l$ rows and $l$ columns indexed by $\mathcal{T}\subset\{1,2,\ldots,n\}$, where $l$ is even and $l=|\mathcal{T}|\le 8$, then for $n\ge n_0$ and $\beta_D/n\le \epsilon_0/50$, $|\mu_B|\le 8.07$, $|\sigma_B^2-1|\le .52$, and $\beta_B\le 50\beta_D$.
\end{lemma}
\begin{proof}
To prove the first claim, we note since $d_{i+}=0, |d_{ij}|\le 1$ and $l\le 8$, letting $\overline{n}=n-l$,
\bea
|b_{i+}|=|\sum_{j\notin\mathcal{T}} d_{ij}|=|-\sum_{j\in\mathcal{T}} d_{ij}|\le 8\quad\mbox{and hence}\quad
\label{bd-b++}
|b_{++}|=|\sum_{i\notin\mathcal{T}} b_{i+}|\le 8\overline{n}.
\ena
This inequality implies $|\mu_B|=|b_{++}|/(\overline{n}-1) \le 8.07$ proving the first claim.

Using $\sigma_D^2=1$ and (\ref{mean-var-Ye}) we obtain
\beq
2\frac{\overline{n}-2}{(\overline{n}-1)(\overline{n}-3)}\sum_{1\le i,j\le n} d^2_{ij}&=& \frac{(\overline{n}-2)(n-1)(n-3)}{(n-2)(\overline{n}-1)(\overline{n}-3)}\in[1,1.01]\quad\mbox{when $n \ge 1000$}\label{nbar}.
\eeq
Using the above bound, (\ref{mean-var-Ye}) and (\ref{bd-b++}), we obtain
\beq
\lefteqn{\left|\sigma^2_B-2\frac{\overline{n}-2}{(\overline{n}-1)(\overline{n}-3)}\sum_{1\le i,j\le n} d^2_{ij}\right|}\\
&\le&\frac{2}{(\overline{n}-1)(\overline{n}-3)}\left((\overline{n}-2)|\sum_{\{i\notin\mathcal{T}\}\cap\{j\notin\mathcal{T}\}} d^2_{ij}-\sum_{i,j=1}^n d^2_{ij}|+\frac{1}{\overline{n}-1}b^2_{++}\right.\nonumber
\left.+2\sum_{i\notin\mathcal{T}} b^2_{i+}\right)\nonumber\\
&\le & \frac{2}{(\overline{n}-1)(\overline{n}-3)}\left((\overline{n}-2)\sum_{\{i\in\mathcal{T}\}\cup\{j\in\mathcal{T}\}} d^2_{ij} + \frac{64\overline{n}^2}{\overline{n}-1}+128\overline{n}\right)\nonumber\\
\label{sum-prime}&\le& \frac{32n}{272(\overline{n}-3)}+\frac{128\overline{n}^2}{(\overline{n}-1)^2(n-3)}+\frac{256\overline{n}}{(\overline{n}-1)(\overline{n}-3)}\\
&\le & .51\quad\mbox{when $n\ge 1000$}\label{sigma-ge-.5},
\eeq

where for (\ref{sum-prime}), we use $|\mathcal{T}|\le 8$, $\beta_D\le n/(50 \times 90)$ and H\"{o}lder's inequality to obtain
\beas
\sum_{\{i\in\mathcal{T}\}\cup\{j\in\mathcal{T}\}}d^2_{ij}&\le& \sum_{i\in\mathcal{T}}\sum_{j=1}^n d^2_{ij}+\sum_{j\in\mathcal{T}}\sum_{i=1}^n d^2_{ij}=2\sum_{i\in\mathcal{T}}\sum_{j=1}^n d^2_{ij}\\
&\le & 2\times 8 \beta_D^\frac{2}{3}n^\frac{1}{3}\quad\mbox{since $\sum_{j=1}^n d^2_{ij}\le (\sum_{j=1}^n |d_{ij}|^3)^\frac{2}{3}n^\frac{1}{3}\le\beta_D^\frac{2}{3}n^\frac{1}{3}$ and $|\mathcal{T}|\le 8$}\\
& \le& 16\frac{n}{272}\quad\mbox{since $\beta_D\le n/4500$}.
\enas
From (\ref{sigma-ge-.5}) and (\ref{nbar}), we obtain $|\sigma^2_B-1|\le .52$ for $n\ge 1000$, proving the second claim of the lemma.

To prove the final claim we observe $|b_{i+}|=|\sum_{j\in\mathcal{T}} d_{ij}|$. Thus,
by (\ref{jensen}) and $|\mathcal{T}|\le 8$,
\bea
&\sum_{i\notin\mathcal{T}} |b_{i+}|^3 = \sum_{i \notin \mathcal{T}} |\sum_{j\in\mathcal{T}} d_{ij}|^3 \le 64 \sum_{i \notin \mathcal{T}}\sum_{j\in\mathcal{T}} |d_{ij}|^3
\le 64\beta_D\quad\mbox{and similarly}\nonumber\\
& |b_{++}|^3=|\sum_{i\notin\mathcal{T}}b_{i+}|^3\le \overline{n}^2\sum_{i\notin\mathcal{T}}|b_{i+}|^3\le 64\overline{n}^2\beta_D.\label{bd-betab-betaa}
\ena
 These observations, along with the fact that $\sigma_B^2\ge .48$, yield
\beas
\beta_B&=&\frac{\sum_{\{i\notin\mathcal{T}\}\cap\{j\notin\mathcal{T}\}} |\widehat{b}_{ij}|^3}{\sigma_B^3}\le  3.1\sum_{\{i\notin\mathcal{T}\}\cap\{j\notin\mathcal{T}\}} |\widehat{b}_{ij}|^3\\
&=&3.1\sum_{\{i\notin\mathcal{T}\}\cap\{j\notin\mathcal{T}\}} \left|b_{ij}-\frac{b_{i+}}{(\overline{n}-2)}-\frac{b_{+j}}{(\overline{n}-2)}+\frac{b_{++}}{(\overline{n}-1)(\overline{n}-2)}\right|^3\\
&\le & 3.1 \times 4^2 \left(\sum |b_{ij}|^3+2\overline{n}\sum_{i\notin\mathcal{T}} \frac{|b_{i+}|^3 }{(\overline{n}-2)^3}+\overline{n}^2\frac{|b_{++}|^3}{((\overline{n}-1)(\overline{n}-2))^3}\right)\quad\mbox{using (\ref{jensen})}\\
&\le & 49.6 \left(\beta_D+128\frac{\overline{n}}{(\overline{n}-2)^3}\beta_D+64\beta_D\frac{\overline{n}^4}{((\overline{n}-1)(\overline{n}-2))^3}\right)\quad\mbox{using (\ref{bd-betab-betaa})}\\
&\le & 50\beta_D\quad\mbox{when $n\ge 1000$, since $l\le 8$}.
\enas
\indent This completes the proof.
\end{proof}
\begin{lemma}\label{lemma-recursion}
Consider a nonnegative sequence $\{a_n, n \in \mathbb{Z}\}$ such that $a_n=0$ for all $n \le 0$ and
\bea \label{an-rec}
a_n\le c+\alpha\max_{l\in\{4,6,8\}} a_{n-l}\quad\mbox{for all $n \ge 1$},
\ena
where $c \ge 0$ and $\alpha\in (0,1/3)$. Then, with $b=\max\{c, a_1(1-3\alpha)\}$, for all $n\ge 1$
$$
a_n\le \frac{b}{1-3\alpha}.
%a_n\le a_1\left(1-\frac{(1-3\alpha)-b}{3\alpha(1-3\alpha)}\right) \quad \mbox{for all $n$.}
$$
\end{lemma}
\begin{proof}
Letting $b_n, n \ge 1$ be given by the recursion
$$
b_{n+1}=3\alpha b_n+b\quad\mbox{for $n \ge 1$, with $b_1=a_1$,}
$$
explicitly solving yields for $n>1$
$$
b_n=c_1(3\alpha)^n+c_2\quad\mbox{where $c_1=\frac{a_1(1-3\alpha)-b}{3\alpha(1-3\alpha)}$ and $c_2=\frac{b}{1-3\alpha}.$}
$$
Since $c_1\le 0$ and $3\alpha<1$, $b_n$ is increasing and bounded above by $c_2$. Hence it suffices to show $a_n\le b_n$.

Since $a_n$ is nonnegative, (\ref{an-rec}) implies that
\bea\label{an-rec-sum}
a_n\le c+\alpha\sum_{l\in\{4,6,8\}} a_{n-l}\quad\mbox{for all $n \ge 1$}.
\ena
We show $a_m \le b_m$ for all $1 \le m \le n$ by induction. When $n \le 4$, since $a_{n-l}=0$ for $l \in \{4,6,8\}$ we have
$$
a_n \le c \le b \le b_n.
$$
Now supposing the claim
is true for some $n \ge 4$, by (\ref{an-rec-sum}) and the monotonicity of $b_n, n \ge 1$, we obtain
\beas
a_{n+1}\le c+\alpha\sum_{l\in\{4,6,8\}} a_{n+1-l}\le b+\alpha\sum_{l\in\{4,6,8\}} b_{n+1-l}\le b+3\alpha b_{n-3}\le b+3\alpha b_n=b_{n+1}.
\enas

This completes the proof.
\end{proof}
\begin{lemma}
\label{lemma:final}
With $\delta^1(\gamma,n)$ defined as in (\ref{def-delta-0}),
\beas
\delta^1(\gamma,n)\le 2804655\frac{\gamma}{n}\quad\mbox{when $n\ge n_0=1000$}.
\enas
\end{lemma}
\begin{proof}
We will consider a smoothed family of indicator functions indexed by $\lambda>0$, namely
\bea
h_{z,\lambda}(x)=\left\{\begin{array}{clcr}
1&\mbox{if $x\le z$}\\
1+(z-x)/\lambda &\mbox{if $x\in(z,z+\lambda]$}\\
0 &\mbox{if $x>z+\lambda$}.
\end{array}
\right.
\ena
Also, define
\beas
h_{z,0}(x)=1_{(-\infty,z]}(x).
\enas
Let $f_{z,\lambda}$ denote the solution to the following Stein equation
\bea\label{stein-eq}
f^\prime(x)-xf(x)=h_{z,\lambda}-\Phi(h_{z,\lambda})\quad\mbox{where $\Phi(h_{z,\lambda})=E(h_{z,\lambda}(Z))$ and $Z\sim\mathcal{N}(0,1)$}.
\ena
We will need the following key inequality about $f_{z,\lambda}$ from \cite{bolt},
\beq
|f^\prime_{z,\lambda}(x+y)-f^\prime_{z,\lambda}(x)|\le |y|\left(1+2|x|+\frac{1}{\lambda}\int_0^1{1_{[z,z+\lambda]}(x+ry)dr}\right)\quad\mbox{for any $\lambda>0$}\label{ineq-bolt}.
\eeq
We consider $D\in M_n^1(\gamma)$ with $\beta_D/n\le \epsilon_0/51$ and $n\ge 1000$. Let $W=\sum_{i=1}^n d_{i\pi(i)}$ as before.
Using the fact that
$$
h_{z,0}\le h_{z,\lambda}\le h_{z+\lambda,0},
$$
and recalling the definition of $\delta_D$ in (\ref{def-delta-a}), we obtain
\bea\label{smoothing}
\delta_D\le \sup_z |E(h_{z,\lambda}(W))-\Phi(h_{z,\lambda})|+\lambda/\sqrt{2\pi}.
\ena
%Taking supremum over $D\in M_n^1(\gamma)$ gives us
%\beq
%\delta^1(\gamma,n)\le \delta^1(\lambda,\gamma,n)+\lambda/\sqrt{2\pi}.\label{smoothing}
%\eeq
 From (\ref{stein-eq}), (\ref{def-zbias}) and $\mbox{Var}(W)=1$,
 we obtain
\beq
\sup_z |E(h_{z,\lambda}(W)-\Phi(h_{z,\lambda})|&=&
\sup_z |E(f_{z,\lambda}'(W))-E(Wf_{z,\lambda}(W))|\nonumber\\
&=&\sup_z |E(f_{z,\lambda}'(W))-E(f_{z,\lambda}'(W^*))|\nonumber\\
&\le &\sup_z E|f_{z,\lambda}'(W^*)-f_{z,\lambda}'(W)|.\label{expln-zbias}
\eeq
 From (\ref{smoothing}), (\ref{expln-zbias}) and (\ref{ineq-bolt}), we obtain the following
\beq
 \delta_D &\le&\sup_z E|f^\prime_{z,\lambda}(W^*)-f^\prime_{z,\lambda}(W)|+\lambda/\sqrt{2\pi}\nonumber\\
 &\le & \sup_z E\left\{|W^*-W|\left(1+2|W|+\frac{1}{\lambda}\int_0^1{1_{[z,z+\lambda]}(W+r(W^*-W))dr}\right)\right\}+\lambda/\sqrt{2\pi}\nonumber\\
  &=&E|W^*-W|+2E(|W||W^*-W|)\nonumber\\
  & &+\sup_z\frac{1}{\lambda}E\left(|W^*-W|\int_0^1{1_{[z,z+\lambda]}(W+r(W^*-W))dr}\right)+\lambda/\sqrt{2\pi}\nonumber\\
 &:=& A_1+A_2+A_3+\lambda/\sqrt{2\pi}\quad\mbox{say}.\label{def-a}
\eeq
First, $A_1$ is bounded using Theorem \ref{thm-w-wstar}, as follows
\beq
A_1&=&E|W^*-W|\le112\frac{\beta_D}{n}+672\frac{\beta_D}{n^2}+192\frac{\beta_D}{n^3}\nonumber\\
&\le &113\frac{\beta_D}{n}\quad\mbox{since $n\ge n_0=1000$}.\label{bd-a1-final}
\eeq
Next we bound $A_2$. From Theorem \ref{thm-w-wstar}, we obtain,
\bea
W^*-W&=&(UW^\dag+(1-U)W^\ddag)-W\nonumber\\
&=&(U(S+T^\dag)+(1-U)(S+T^\ddag))-(S+T)=UT^\dag+(1-U)T^\ddag-T.\label{diff-star}
\ena
Let $\mathbf{I}=(I^\dag,J^\dag,K^\dag,L^\dag,\pi(I^\dag),\pi(J^\dag),\pi(K^\dag),\pi(L^\dag))$ and $\mathcal{I}$ be as defined in Lemma \ref{lemma-pi-dag-involution}. Thus, by (\ref{composition-def}), the right hand side of (\ref{diff-star}), and hence $W^*-W$, is measurable with respect to $\mathscr{I}=\{\mathbf{I},U\}$. Furthermore, since $\sup_{i,j}|d_{ij}|\le 1$ and $|\mathcal{I}|\le 8$, we have $$|W|=|S+T|\le |S|+|T|=|S|+|\sum_{i\in\mathcal{I}}d_{i\pi(i)}|\le |S|+\sum_{i\in\mathcal{I}}|d_{i\pi(i)}|\le |S|+8.$$
Now, by the definition of $A_2$, and that $U$ is independent of $S$ and ${\bf I}$,
\beq
A_2&=& 2E\left(|W||W^*-W|\right)\nonumber\\
&=& 2E\left(|W^*-W|E(|W||\mathscr{I})\right)\le 2E\left(|W^*-W|E(|S|+8|\mathscr{I})\right)\nonumber\\
&\le & 2E\left(|W^*-W|\sqrt{E(S^2|\mathbf{I})}\right)+16E|W^*-W|.\label{bd-a2-1}
\eeq
In the following, for $\mathbf{i}$ a realization of $\mathbf{I}$, let $|\mathbf{i}|$ denote the number of
distinct elements of $\mathbf{i}$. Since $S=\sum_{i\notin\mathcal{I}} d_{i\pi(i)}$ and $\pi$ is chosen from $\Pi_n$ uniformly at random, conditioned on $\mathbf{I}=\mathbf{i}$, $S$ has the same distribution as $\sum_{i\notin \mathbf{i}} b_{i\theta(i)}$, where $B=((b_{ij}))_{n-|\mathbf{i}|\times n-|\mathbf{i}|}$ is the matrix formed by removing rows and columns  of $D$ that occur in $\mathbf{i}$ and $\theta$ is chosen uniformly from $\Pi_{n-|\mathbf{i}|}$. Since $|\mathbf{i}|\le 8$ and $\beta_D/n\le\epsilon_0/51<\epsilon_0/50$, Lemma \ref{lemma-smaller-b} yields $|\mu_B| \le 8.07$ and $\sigma^2_B \le 1.52$, and hence
\beq
E|Y_B|^2\le 1.52+8.07^2=66.6449\quad\mbox{when $n\ge 1000$}.\label{bd-var-S}
\eeq
Using (\ref{bd-var-S}), we obtain
\beas
E(S^2|\mathbf{I}=\mathbf{i})= E|Y_B|^2\le 67 \quad \mbox{for all $\mathbf{i}$.}
\enas
Thus using (\ref{bd-a2-1}) and (\ref{bd-a1-final}), we obtain
\beq
A_2&\le & 33 E|W^*-W|\le 3729\frac{\beta_D}{n}.\label{bd-a2-final}
\eeq

Finally, we are left with bounding $A_3$. First we note that
\beas
W+r(W^*-W)&=&rW^*+(1-r)W\\
&=&r(S+UT^\dag+(1-U)T^\ddag)+(1-r)(S+T)\\
&=& S+rUT^\dag+r(1-U)T^\ddag+(1-r)T\\
&=& S+g_r\quad\mbox{where $g_r=rUT^\dag+r(1-U)T^\ddag+(1-r)T$}.
\enas
Now, from the definition of $A_3$, again using that $W-W^*$ is $\mathscr{I}$ measurable,
\beq
A_3&=&\sup_z \frac{1}{\lambda}E\left(|W-W^*|\int_0^1{1_{[z,z+\lambda]}(W+r(W^*-W))dr}\right)\nonumber\\
&=&\sup_z \frac{1}{\lambda}E\left(|W-W^*|E(\int_0^1{1_{[z,z+\lambda]}(W+r(W^*-W))dr|\mathscr{I}})\right)\nonumber\\
&=&\sup_z\frac{1}{\lambda}E\left(|W-W^*|\int_0^1P(W+r(W^*-W)\in[z,z+\lambda]|\mathscr{I})dr\right)\nonumber\\
&=&\sup_z\frac{1}{\lambda}E(|W-W^*|\int_0^1P(S+g_r\in[z,z+\lambda]|\mathscr{I})dr)\nonumber\\
&=&\sup_z\frac{1}{\lambda}E(|W-W^*|\int_0^1P(S\in[z-g_r,z+\lambda-g_r]|\mathscr{I})dr)\nonumber\\
\nonumber &\le&\frac{1}{\lambda}E(|W-W^*|\int_0^1\sup_z P(S\in[z-g_r,z+\lambda-g_r]|\mathscr{I})dr)\\
  &=&\frac{1}{\lambda}E(|W-W^*|\int_0^1\sup_z P(S\in[z,z+\lambda]|\mathscr{I})dr)\label{mbility}\\
&=&\frac{1}{\lambda}E(|W-W^*|\sup_zP(S\in[z,z+\lambda]|\mathscr{I}))\nonumber \\
&=&\frac{1}{\lambda}E(|W-W^*|\sup_z P(S\in[z,z+\lambda]|\mathbf{I})),\label{bd-a3}
\eeq
where to obtain equality in (\ref{mbility}) we have used the fact that $g_r$ is measurable with respect to $\mathscr{I}$ for all $r$.

It remains only to bound $P(S\in[z,z+\lambda]|\mathbf{I})$. In the following calculations $\widetilde{b}_{ij}={\widehat{b}_{ij}}/{\sigma_B}$ as before.
Since $\beta_D/n\le \epsilon_0/51$, Lemma \ref{lemma-smaller-b} yields $\sigma_B>1/2$. Hence,
\beq
\nonumber \sup_z P(S\in[z,z+\lambda]|\mathbf{I}=\mathbf{i}))&=&\sup_zP(\sum_{i\notin\mathbf{i}} b_{i\theta(i)}\in[z,z+\lambda])\\
&=&\sup_zP(\sum_{i\notin\mathbf{i}}\frac{b_{i\theta(i)}}{\sigma_B}\in [\frac{z}{\sigma_B},\frac{z+\lambda}{\sigma_B}])\nonumber\\
&\le &\sup_z P(\sum_{i\notin\mathbf{i}} \frac{b_{i\theta(i)}}{\sigma_B}\in[\frac{z}{\sigma_B},\frac{z}{\sigma_B}+2\lambda])\nonumber\\
&=&\sup_z P(\sum_{i\notin\mathbf{i}} \frac{b_{i\theta(i)}}{\sigma_B}\in [z,z+2\lambda])\nonumber\\
&=& \sup_z P(\sum_{i\notin\mathbf{i}} \widetilde{b}_{i\theta(i)}\in[z,z+2\lambda])\nonumber \\
&=& \sup_z P(Y_{\widetilde{B}} \in[z,z+2\lambda])\label{centering}.
\ena
The equality (\ref{centering}) holds as when computing ${\widehat{b}_{ij}}$ we have that $\sum_{i\notin\mathbf{i}}b_{i+}$ and $\sum_{i\notin\mathbf{i}} b_{+\theta(i)}=\sum_{j\notin\mathbf{i}} b_{+j}$ do not depend on $\theta.$ Recalling that the distribution function of $Y_{\widetilde{B}}$ is denoted by $F_{\widetilde{B}}$, we have, from the definition of $\delta(\gamma,n)$ in (\ref{def-delta}),
\bea
P(Y_{\widetilde{B}}\in[z,z+2\lambda])&\le &|F_{\widetilde{B}}(z+2\lambda)-\Phi(z+2\lambda)|+|F_{\widetilde{B}}(z)-\Phi(z)|+|\Phi(z+2\lambda)-\Phi(z)|\nonumber\\
&\le & 2\delta_{\widetilde{B}}+\frac{2\lambda}{\sqrt{2\pi}}.\label{ineq-last}
\ena

Note that $\beta_{\widetilde{B}}=\beta_B$ and by Lemma \ref{lemma-smaller-b}
$$
\frac{\beta_{\widetilde{B}}}{n-|\mathbf{i}|}=\frac{\beta_B}{n-|\mathbf{i}|}\le 50\frac{\beta_D}{n-|\mathbf{i}|}\le 51\frac{\beta_D}{n}\le \epsilon_0.
$$
Hence using the fact that $\widetilde{B}\in M_{n-|\mathbf{i}|}(\beta_B)$ and applying Lemma \ref{bd-del-del1}, we obtain for $n\ge 1008$
\bea\label{ineq-last-1}
\delta_{\widetilde{B}}\le \delta^1(22\beta_{\widetilde{B}})+36\frac{\beta_{\widetilde{B}}}{{n-|\mathbf{i}|}}=\delta^1(22\beta_B,n-|\mathbf{i}| )+36\frac{\beta_{{B}}}{{n-|\mathbf{i}|}}.
\ena
Inequalities (\ref{ineq-last}) and (\ref{ineq-last-1}) imply
\bea
P(Y_{\widetilde{B}}\in[z,z+2\lambda])&\le & 2\delta_{\widetilde{B}}+\frac{2\lambda}{\sqrt{2\pi}}\nonumber\\
&\le& 2\delta^1(22\beta_B,n-|\mathbf{i}|)+72\frac{\beta_B}{n-|\mathbf{i}|}+\frac{2\lambda}{\sqrt{2\pi}}\nonumber\\
&\le& 2\delta^1(22\times 50\beta_D,n-|\mathbf{i}|)+\frac{72\times 50\beta_D}{n-|\mathbf{i}|}+\frac{2\lambda}{\sqrt{2\pi}}\nonumber\\
&\le& 2\max_{l\in\{4,6,8\}}\delta^1(k_1\beta_D,n-l)+\frac{k_2\beta_D}{n}+\frac{2\lambda}{\sqrt{2\pi}}\label{concentration},
\eeq
where $k_1=1100$ and $k_2=3630$. As (\ref{concentration}) does not depend on $z$ or $\mathbf{i}$, it bounds $\sup_z P(S\in[z,z+\lambda]|\mathbf{I}))$ in (\ref{bd-a3}).
\par Now using (\ref{bd-a3}), (\ref{concentration}) and (\ref{bd-a1-final}), we obtain
\beq
A_3&\le &\frac{1}{\lambda}(2\max_{l\in\{4,6,8\}}\delta^1(k_1\beta_D,n-l)+k_2\frac{\beta_D}{n}+\frac{2\lambda}{\sqrt{2\pi}})E|W-W^*|\nonumber\\
&\le &\frac{1}{\lambda}(2\max_{l\in\{4,6,8\}}\delta^1(k_1\beta_D,n-l)+k_2\frac{\beta_D}{n}+\frac{2\lambda}{\sqrt{2\pi}})113\frac{\beta_D}{n}\label{bd-a3-final}.
\eeq
Combining (\ref{bd-a1-final}), (\ref{bd-a2-final}), (\ref{bd-a3-final}) and using (\ref{def-a}) we obtain for $n\ge 1008$,
\beas
\delta_D&\le& 3842\frac{\beta_D}{n}+\frac{1}{\lambda}(2\max_{l\in\{4,6,8\}}\delta^1(k_1\beta_D,n-l)+k_2\frac{\beta_D}{n}+\frac{2\lambda}{\sqrt{2\pi}})113\frac{\beta_D}{n}+\frac{\lambda}{\sqrt{2\pi}}\\
&\le& 3842\frac{\gamma}{n}+\frac{113\gamma}{n\lambda}(2\max_{l\in\{4,6,8\}}\delta^1(k_1\beta_D,n-l)+k_2\frac{\gamma}{n}+\frac{2\lambda}{\sqrt{2\pi}})+\frac{\lambda}{\sqrt{2\pi}},
\enas
since $\beta_D \le \gamma$ for all $D \in M_n^1(\gamma)$.
Setting $\lambda=(113 \times 8) k_1\gamma/n$, we obtain for $n\ge 1008$,
\bea\label{bd-d-first}
\delta_D&\le & c\frac{\gamma}{n}+\frac{1}{4}\max_{l\in\{4,6,8\}}\frac{\delta^1(k_1\beta_D,n-l)}{k_1},
\ena
where $c=400,665$. Since $c>51/\epsilon_0$, its clear that if we consider $D\in M_n^1(\gamma)$, with $\beta_D/n>\epsilon_0/51$, then $\delta_D< c\gamma/n$. Hence (\ref{bd-d-first}) holds for all $D\in M_n^1(\gamma)$ with $n\ge 1008$. Taking supremum over $D\in M_n^1(\gamma)$, we have, for $n \ge 1008$,
\beas
\delta^1(\gamma,n)& \le & c\frac{\gamma}{n}+\frac{1}{4}\max_{l\in\{4,6,8\}}\frac{\delta^1(k_1\gamma,n-l)}{k_1},
\enas
where $c=400,665$. Now multiplying by $n/\gamma$ and taking supremum over $\gamma$ we obtain
\beq
\sup_\gamma n\frac{\delta^1(\gamma,n)}{\gamma}\le c+\frac{2}{7}\max_{l\in\{4,6,8\}}\sup_\gamma (n-l)\frac{\delta^1(\gamma,n-l)}{\gamma}\quad\mbox{for all $n\ge 1008$}.\label{eqn-recursion-8}
\eeq
If $D\in M^1_{n}(\gamma)$ and $n\ge 1000$ then (\ref{bd-betad}) shows that $\beta_D  \ge 2$, and hence $D \not \in M^1_{n}(\gamma)$ for all
$\gamma < 2$.
Since $\delta^1(\gamma,n)\le 1$ for all $n \in \mathbb{N}$, for $1000\le n\le 1008$ we have
$$
\sup_\gamma  n\frac{\delta^1(\gamma,n)}{\gamma}=\sup_{\gamma \ge 2} n\frac{\delta^1(\gamma,n)}{\gamma}\le 1008\sup_{\gamma \ge 2}\frac{\delta^1(\gamma,n)}{\gamma}\le504.
$$
Since $c>504$, we conclude (\ref{eqn-recursion-8}) holds for $n\ge 1000$, that is
\beq
\sup_\gamma n\frac{\delta^1(\gamma,n)}{\gamma}\le c+\frac{2}{7}\max_{l\in\{4,6,8\}}\sup_\gamma (n-l)\frac{\delta^1(\gamma,n-l)}{\gamma}\quad\mbox{for all $n\ge 1000$}.\label{eqn-recursion}
\eeq

Letting $s_n=\sup_\gamma n\delta^1(\gamma,n)/\gamma$ and $a_n=s_{n+999}$ for $n\ge 1$, and $a_n=0$ for $n\le 0$, (\ref{eqn-recursion}) gives
\beas
a_n\le c +\frac{2}{7}\max_{l\in\{4,6,8\}} a_{n-l}\quad\mbox{for $n\ge 1$.}
\enas
%where for $1 \le n \le 8$ the inequality holds again using $c>504$.
Using Lemma \ref{lemma-recursion} with $\alpha=2/7$, $c=400,665$ and
noting that $a_1=n\delta^1(\gamma,1000)/\gamma\le 500$ since $\gamma \ge 2$, we obtain
$$
a_1(1-3\alpha) \le 500(1-6/7) < c \quad \mbox{which yields} \quad b=\max\{c,a_1(1-3\alpha)\}=c,
$$
and therefore
$$
a_n \le \frac{b}{1-3\alpha}=\frac{c}{1-3\alpha}= 2804655.
$$
This completes the proof of the lemma.
\end{proof}

Lemmas \ref{lemma:final}, \ref{bd-del-del1} and remarks following (\ref{bd-betad}) yield
\beas
||F_W-\Phi||_\infty\le K\frac{\beta_D}{n}\quad\mbox{with $n\ge 1000$},
\enas
where we can take $K=22\times 2804655+36=61,702,446$. This proves Theorem \ref{prop-main} and hence completes the proof of Theorem \ref{Lp-theorem} as well.
\section*{Appendix}
Here we briefly indicate that the order $\beta_D/n$ in Theorem \ref{prop-main} can not be improved uniformly over all arrays $D$ satisfying the conditions as in the theorem. For example, define the symmetric array $E$ given as follows where because of symmetry we define the entries $e_{ij}$ for $j\ge i$ only.
\beas
e_{ij}=\left\{\begin{array}{ll}
0 & \mbox{if $i=j$ or $i$ is odd and $j=i+1$}\\
1 & \mbox{if $i\neq j$ and $i-j$ is even}\\
-1 & \mbox{otherwise}.
\end{array}.
\right.
\enas
Clearly $e_{i+}=0$ for all $1\le i\le n$ and for $i=2k-1$
\beas
\sum_{j\ge i+1} e^2_{i+1 j}=\sum_{j\ge i} e^2_{ij}=n-2k.
\enas
Using symmetry again, we obtain $\sum_{1\le i,j\le n} e^2_{ij}=2\sum_{i, j\ge i} e^2_{ij}=O(n^2)$ and using (\ref{sigmaD}), we have $\sigma^2_E=O(n)$. Also, since $|e_{ij}|=0,1$, we have $\beta_E=f_n/g_n$, where $f_n=\sum |e^3_{ij}|=O(n^2)$ and $g_n=\sigma^3_E=O(n^{3/2})$. Collecting all these facts together, we obtain $\beta_E/n=O(n^{-1/2})$. Also, define $D$ by $d_{ij}=e_{ij}/\sigma_E$.

Now along the lines of \cite{golpen}, fix $\epsilon\in (0,1)$ and define $t=(1-\epsilon)/\sigma_E$. Then, we have
\beas
\Phi(t)-\Phi(0)\ge t \phi(t)\ge \left(\frac{1-\epsilon}{\sigma_E}\right)\phi(\sigma_E^{-1}).
\enas
Using notations as in (\ref{def-U}) and (\ref{def-Y}), we observe that $Y_E$ is integer valued implying $F_W(0)=F_W(t)$ and hence
\beas
||F_W-\Phi||_\infty\ge \frac{1}{2}(1-\epsilon)\sigma_E^{-1}\phi(\sigma_E^{-1}).
\enas
Multiplying by $n^{1/2}$ on both sides yields
\bea
n^{1/2}||F_W-\Phi||_\infty\ge n^{1/2}\frac{1}{2}(1-\epsilon)\sigma_E^{-1}\phi(\sigma_E^{-1}).\label{laststep}
\ena
Since $\sigma_E=O(n^{1/2})$, letting $n\rightarrow\infty$ and then taking $\epsilon\rightarrow0$, we obtain
\beas
\liminf_{n\rightarrow\infty} n^{1/2}||F_W-\Phi||_\infty\ge O(1).
\enas
Since $\beta_E/n=O(n^{-1/2})$, we have that $\beta_E/n$ provides the correct rate of convergence.
\section*{Acknowledgements}
The author would like to thank his advisor Prof. Larry Goldstein for suggesting the problem, pointing out several mistakes in earlier working versions of this paper and giving valuable insight.

 \end{document}